# Weak-form modified sparse identification of nonlinear dynamics

Cristian López[a], Ángel Naranjo[a], Diego Salazar[a], Keegan J. Moore[a,b,]*

[a]Department of Mechanical and Materials Engineering, University of Nebraska, Lincoln, NE 68588, United States
[b]Daniel Guggenheim School of Aerospace Engineering, Georgia Institute of Technology, Atlanta, GA 30332

* Corresponding author.

e-mail addresses: clopez19@huskers.unl.edu (C. López); anaranjorubio2@huskers.unl.edu (Á. Naranjo); dsalazarorellana2@huskers.unl.edu (D. Salazar); kmoore@gatech.edu (K.J. Moore).

## Abstract

Identifying nonlinear dynamics and characterizing noise from data is critical across science and engineering for understanding and modeling the behavior of the systems accurately. The modified sparse identification of nonlinear dynamics (mSINDy) has emerged as an effective framework for identifying systems embedded in heavy noise; however, further improvements can expand its capabilities and robustness. By integrating the weak SINDy (WSINDy) into mSINDy, we introduce the weak mSINDy (WmSINDy) to improve the system identification and noise modeling by harnessing the strengths of both approaches. The proposed algorithm simultaneously identifies parsimonious nonlinear dynamics and extracts noise probability distributions using automatic differentiation. We evaluate WmSINDy using several nonlinear systems and it demonstrates improved accuracy and noise characterization over baselines for systems embedded in relatively strong noise.

**Notation**

| Symbol | Meaning | Domain |
|---|---|---|
| $(\mathbf{B}, \mathbf{G})$ | WSINDy linear system | $\mathbb{R}^{H \times D}, \mathbb{R}^{H \times J}$ |
| $D$ | dimension | $\mathbb{Z}^+$ |
| $\boldsymbol{f}$ | nonlinear function | $\mathbb{R}^D$ |
| $\mathcal{F}$ | Fourier transform | $\mathbb{R}^D$ |
| $\widetilde{\mathbf{F}}$ | estimated states | $\mathbb{R}^{N \times D}$ |
| $k^*$ | wavenumber | $\mathbb{Z}^+$ |
| $J$ | number of trial functions in library | $\mathbb{Z}^+$ |
| $\lambda$ | threshold for sparsity promotion | $\mathbb{R}^+$ |
| $\mathcal{L}_{mSINDy}\left(\boldsymbol{\Xi}, \widetilde{\mathbf{N}}(t)\right)$ | mSINDy loss function | $BV(0, \infty)$ |
| $\mathcal{L}\left(\boldsymbol{\Xi}, \widetilde{\mathbf{N}}(t)\right)$ | WmSINDy loss function | $BV(0, \infty)$ |
| $m$ | support of $\phi$ | $\mathbb{Z}^+$ |
| $N$ | number of points | $\mathbb{Z}^+$ |
| $N_{loops}$ | iteration loops | $\mathbb{Z}^+$ |
| $\mathbf{N}(t)$ | measurement noise | $\mathbb{R}^{N \times D}$ |
| $\widetilde{\mathbf{N}}(t)$ | learned noise | $\mathbb{R}^{N \times D}$ |
| $p$ | degree of $\phi$ | $\mathbb{Z}^+$ |
| $\phi$ | test function of compact support in the time domain | $C^1(\mathbb{R}^D)$ |
| $\hat{\phi}$ | test function in the frequency domain | $C(\mathbb{R}^D)$ |
| $q$ | prediction step | $\mathbb{Z}^+$ |
| $\boldsymbol{\Theta}$ | library of functions | $\mathbb{R}^{N \times J}$ |
| $t_T$ | total time | $\mathbb{R}^+$ |
| $t$ | time | $\mathbb{R}^N$ |
| $\tau$ | decay tolerance of $\phi$ in the real space | $\mathbb{R}^+$ |
| $\hat{\tau}$ | decay tolerance of $\hat{\phi}$ in the Fourier space | $\mathbb{R}^-$ |
| $\mathbf{x}(t)$ | system states | $\mathbb{R}^D$ |
| $\mathbf{X}(t)$ | true solution | $\mathbb{R}^{N \times D}$ |
| $\widetilde{\mathbf{X}}(t)$ | estimated solution | $\mathbb{R}^{N \times D}$ |
| $\boldsymbol{\Xi}$ | true sparse matrix of coefficients | $\mathbb{R}^{J \times D}$ |
| $\widetilde{\boldsymbol{\Xi}}$ | learned sparse matrix of coefficients | $\mathbb{R}^{J \times D}$ |
| $\mathbf{Y}(t)$ | measurement data | $\mathbb{R}^{N \times D}$ |

# 1. Introduction

Understanding the fundamental physics that governs intricate systems is a pivotal task in various fields of science and engineering to model the behavior of these systems accurately[1]. In traditional physics and engineering, governing equations are derived from first principles such as conservation of energy and momentum. Nonetheless, in many modern systems, the derivation of governing equations from these principles may be quite challenging or even impossible [2]. This could be due to the complexity of the system, incomplete comprehension of the underlying mechanics, or the absence of universal laws. Currently, advances in big data collection, sensing technologies, machine learning, etc. have made it possible to model complex dynamical systems from measurements [3]. In this sense, discovering the governing equations underlying a dynamical system simply from data measurements provides a promising path to gain insights into the behavior of a system [4], make reliable predictions [5], and develop reduced-order models [6][7] that are interpretable and generalizable.

      The sparse identification of nonlinear dynamics (SINDy), proposed by Brunton et al. [2], is a significant milestone in identifying meaningful mathematical models that accurately describe the observed behavior of a system. After collecting the sample measurements, the analyst provides a library of possible functions that could appear in the system's equations, then it performs a sequentially thresholded, least-squares regression to recover the fewest terms required to represent the data. The regression is performed between the numerical derivative of the data and the library of functions. SINDy has been applied in various fields, such as in epidemiology to investigate the transmission dynamics of COVID-19 for constructing a global model [8] in biology for representing biological networks [9], and in biophysical process modeling for learning stochastic equations [10]. Other examples include processes modeling and control in chemical engineering [11], correcting nonlinear impairments in fiber optic transmissions [12], modeling low-dimensionalized complex fluid flows [13], and identifying nonlinear and hysteric behaviors [14]. The algorithm provides accurate results in both noise-free environments and slight presence of noise, however, produces complex solutions when the data is embedded in strong noise. In this regard, several improvements have been proposed to extend its application. Namely, to construct the correct library of terms [14–18], to improve the sparse regression [19–25], and to deal with strong noisy data [5,17,26,27].

      Although successful, the SINDy algorithm [2] and its improvements [5,14,16-28] rely on numerical derivatives, which limits the performance when the data has a low signal-to-noise ratio. There have been several research efforts put forth to address this issue, which replace the derivatives with integrals in the

sparse model selection [28–32] to form a weak form of SINDy. As explained in [30,31], the weak formulation applied to SINDy has been shown to provide a more robust way to handle noisy data, because the derivative acts on smooth and localized test functions. This acts as a natural filtering/regularization that reduces the effect of noise in the data. Thus, it avoids computing numerical derivatives of noisy data that can amplify errors. In [33], it was used a neural network with a Runge-Kutta integration framework to constrain the learning algorithm inside a numerical time-stepping scheme. The resulting method simultaneously learns the dynamics and estimates the measurement noise by treating both the measurement error and the dynamics as unknowns to be identified. Following this research, automatic differentiation [34] and time-stepping constraints were combined in a SINDy framework to simultaneously denoise data, learn and parametrize noise probability distributions, and identify underlying parsimonious dynamical systems responsible for generating the time-series data [26]. This algorithm coined modified SINDy (mSINDy), handles significant noise levels in the data, sometimes up to 40% of the amplitude of the embedded signal. However, mSINDy still applies numerical derivatives to the data and could be improved further.

In this study, we introduce the WmSINDy, which leverages the inherent properties of the WSINDy and mSINDy methods. The WmSINDy method is built on the mSINDy framework to improve the identification of system dynamics using signals with heavy noise, an endeavor previously challenging for the conventional mSINDy approach. The significance of the proposed approach is that it offers greater insight into the nature of noise in dynamical systems and how it affects equation discovery, yields more accurate identification of equation coefficients, and adapts to various noise distributions. Thus, the proposed method is applicable to a wide range of real-world scenarios. We showcase the approach on several numerical examples from canonical systems such as the Duffing oscillator to systems that can exhibit chaos. The paper is organized as follows. In Section 2, we present the methods used. Section 3 introduces the proposed framework. In Section 4, we demonstrate performance on benchmark systems. Section 5 provides the conclusion and open questions.

## 2. Methods

### 2.1 SINDy

Consider a first-order dynamical system in $D$ dimensions of the form

$$\frac{d}{dt}\mathbf{x}(t) = \boldsymbol{f}(\mathbf{x}(t)), \qquad \mathbf{x}(0) = \mathbf{x}_0 \in \mathbb{R}^D, \qquad 0 \leq t \leq t_T, \tag{1}$$

where $\mathbf{x}(t) = [x_1(t)\ x_2(t)\ \ldots\ x_D(t)]^\mathrm{T} \in \mathbb{R}^D$ is the system state, $\boldsymbol{f}$ is in general a nonlinear function that governs the evolution of $\mathbf{x}(t)$, and $t_T$ the total time. The SINDy algorithm is successful in solving this problem for sparsely represented nonlinear dynamics when noise is small, and the dynamic scales do not vary across multiple orders of magnitude. SINDy asserts the ansatz that a set of candidate functions characterize the right-hand side of the governing equations by constructing a library of functions $\theta_j$, $j = 1,2,\ldots,J$, such that $\boldsymbol{\Theta}(\mathbf{X}) \in \mathbb{R}^{N \times J}$ is

$$\boldsymbol{\Theta}(\mathbf{X}) = [\theta_1(\mathbf{X})\ \theta_2(\mathbf{X}) \ldots \theta_J(\mathbf{X})],$$

where $\mathbf{X} = [\mathbf{x}(t_1)\ \mathbf{x}(t_2)\ \ldots\ \mathbf{x}(t_N)]^\mathrm{T} \in \mathbb{R}^{N \times D}$, and $\theta_j(\mathbf{X})$ is a point-wise evaluation. Thus, to select the few important basis functions to describe the dynamics, a sparse regression problem is constructed as

$$\Xi = \arg\min_{\Xi} \left\| \dot{\mathbf{X}} - \boldsymbol{\Theta}(\mathbf{X})\Xi \right\|_2 + \lambda \|\Xi\|_0, \tag{2}$$

where $\Xi = [\xi_1\ \xi_2\ ...\ \xi_d] \in \mathbb{R}^{J \times D}$ is a sparse weight matrix obtained by using sequential-thresholding least squares controlled by the $\lambda > 0$ hyperparameter. By solving Eq. (2), the model of Eq. (1) is expressed as

$$\frac{d}{dt}\mathbf{x}(t) = f(\mathbf{x}(t)) \approx \Theta(\mathbf{x}(t))\Xi. \tag{3}$$

**2.2 Modified SINDy**

By using automatic differentiation [34], the mSINDy has two particularities that make it an important advancement for the data-driven discovery of dynamical systems: first, the algorithm simultaneously identifies the parsimonious equations that govern the dynamical system; and second, it characterizes the noise $\widetilde{\mathbf{N}}$ that can be separated from the signals. These features make mSINDy better at identifying dynamical systems from noisy measurements than the standard SINDy algorithm.

In real scenarios, the measurement data $\mathbf{Y} = [\mathbf{y}(t_1)\ \mathbf{y}(t_2)\ ...\ \mathbf{y}(t_N)]^T \in \mathbb{R}^{N \times D}$, contains noise such that $\mathbf{Y}(t) = \mathbf{X}(t) + \mathbf{N}(t)$, where $\mathbf{N} = [\mathbf{n}(t_1)\ \mathbf{n}(t_2)\ ...\ \mathbf{n}(t_N)]^T \in \mathbb{R}^{N \times D}$ is a matrix of the noise, which makes the solution of Eq. (2) more challenging. In this sense, by leveraging the automatic differentiation of Tensorflow [35] to simultaneously denoise and learn a model, mSINDy solves the following equation

$$\mathcal{L}_{mSINDy}(\Xi, \mathbf{N}) = e_d + e_s, \tag{4a}$$

$$\mathcal{L}_{mSINDy}(\Xi, \widetilde{\mathbf{N}}) = \left\|\dot{\widetilde{\mathbf{X}}} - \Theta(\widetilde{\mathbf{X}})\Xi\right\|_2^2 + \sum_{j=q+1}^{N-q}\sum_{i=-q, i\neq 0}^{q} \omega_i \left\|\mathbf{y}_{j+i} - \widetilde{\mathbf{n}}_{j+i} - \widetilde{\mathbf{F}}^i(\widetilde{\mathbf{x}}_j)\right\|_2^2, \tag{4b}$$

where $e_d$ is the derivative error, $e_s$ is the simulated error, and $\omega_i = 0.9^{|i|-1}$ ($i \in \mathbb{Z}^+$ and $i \neq 0$ to ensure that the interaction between point $j$ and itself is not doubly counted, preventing redundancy) is a hyperparameter utilized to address numerical errors that may arise during the simulation of the dynamical system. The hyperparameter $\omega_i$ is implemented to ensure that errors at later times are less significant than those at early times in the simulations. In machine learning, the $\|\cdot\|_2^2$ is prevalent because it is continuously differentiable, which facilitates optimization processes[36]. The estimated state $\tilde{\mathbf{x}}$ is obtained by integrating Eq. (3), such as with the 4[th]-order Runge-Kutta method, forward/backward in time $q$-steps as

$$\mathbf{y}(j+q) - \widetilde{\mathbf{n}}(j+q) = \tilde{\mathbf{x}}(j+q) = \widetilde{\mathbf{F}}^q(\tilde{\mathbf{x}}(j)) = \tilde{\mathbf{x}}(j) + \int_{t_j}^{t_{j+q}} \Theta(\tilde{\mathbf{x}}(\tau))\Xi\, d\tau. \tag{5}$$

The optimization problem of Eq. (4), searches for the solution of the coupled parameters $\Xi$ and $\widetilde{\mathbf{N}}$ as

$$\Xi, \widetilde{\mathbf{N}} = \arg\min_{\Xi, \widetilde{\mathbf{N}}} \mathcal{L}(\Xi, \widetilde{\mathbf{N}}),$$

$$\text{s.t.} (|\Xi| < \lambda) = 0. \tag{6}$$

**2.3 Weak SINDy**

The WSINDy is based on the weak formulation of the dynamical system that enables the use of a circular convolutional weak form, such that the Fourier theorem applies. Thus, WSINDy benefits from being able to use the discrete Fourier transform (DFT) and, more specifically, the Fast Fourier Transform (FFT). The DFT is applicable here because the sampling of the dynamics results in discrete data that can be analyzed in the frequency domain. The circular convolutional weak form assumes periodicity in the data, which aligns with the DFT's assumption of periodicity in its input. Thus, the use of FFT not only permits efficient handling of these convolutions but also ensures computational speed, especially when working with large datasets. Through the combination of these, the WSINDy is a powerful and efficient approach for the identification of nonlinear dynamical systems from measured data. The weak formulation is obtained by multiplying (1) by a suitable function $\phi$

$$\int_a^b \frac{d}{dt}\mathbf{x}(t)\,\phi(t)dt = \int_a^b \boldsymbol{f}(\mathbf{x}(t))\,\phi(t)dt. \tag{7}$$

Using Green's theorem,

$$\phi(t_b)\mathbf{x}(t_b) - \phi(t_a)\mathbf{x}(t_a) = \int_a^b \mathbf{x}(t)\phi'(t)dt + \int_a^b \boldsymbol{f}(\mathbf{x}(t))\,\phi(t)dt. \tag{8}$$

Taking $\phi$ to be compactly supported in $(a,b) \subset [0, t_T]$[30]

$$0 = \int_a^b [\phi'(t)\mathbf{x}(t) + \phi(t)\boldsymbol{f}(\mathbf{x}(t))]dt. \tag{9}$$

Note that, for each component of $\mathbf{x}(t)$, there is a corresponding family of test functions. Additionally, restricting the test functions family to translations of a fixed function $\phi$, $\phi_n(t) = \phi(t - t_n)$ for each component, where $t_n = n\Delta t \in [0, t_T]$, with $\Delta t = t_T/(N-1)$, then,

$$0 = \int_a^b [\phi'_n(t)\mathbf{x}(t) + \phi_n(t)\boldsymbol{f}(\mathbf{x}(t))]dt. \tag{10}$$

The right-hand side of Eq. 10 can be interpreted as the sum of the discrete convolution operations: one between $\phi'$ and $\mathbf{x}$, and the other between $\phi$ and $\boldsymbol{f}$. Replacing $\mathbf{x}(t)$ for $\mathbf{Y}(t)$ and $\boldsymbol{f}(\mathbf{x}(t))$ for $\boldsymbol{\Theta}(\mathbf{Y}(t))\boldsymbol{\Xi}$, in a compact notation, a general residual $\mathcal{R}(\boldsymbol{\Xi}; \phi)$ is defined as [30]

$$\mathcal{R}(\boldsymbol{\Xi}; \phi) = \int_a^b [\phi'_n(t)\mathbf{Y}(t) + \phi_n(t)\boldsymbol{\Theta}(\mathbf{Y}(t))\boldsymbol{\Xi}]dt. \tag{11}$$

Then, we define

$$\mathbf{B}(\mathbf{Y}(t)) = -\langle \phi'(t), \mathbf{Y}(t) \rangle \text{ and } \mathbf{G}(\mathbf{Y}(t)) = \langle \phi(t), \boldsymbol{\Theta}(\mathbf{Y}(t)) \rangle \in \mathbb{R}^{H \times J},$$

where $\langle \cdot, \cdot \rangle$ stands for the corresponding convolution operation for each component in the data set, and $H$ is the dimension of the product of the convolution for each component. By using these numerical integrals of the data instead of derivatives, the influence of noise is reduced. Thus, to obtain the governing dynamics of the system, the linear system of equations of the discrete version of the weak form is

$$\mathbf{B}(\mathbf{Y}(t)) = \mathbf{G}(\mathbf{Y}(t))\boldsymbol{\Xi}. \tag{12}$$

The solution of Eq. 12 minimizes $\mathcal{R}$. Thus, to find the coefficients of $\boldsymbol{\Xi} \in \mathbb{R}^{J \times D}$, the following sparse regression problem

$$\boldsymbol{\Xi} = \arg\min_{\boldsymbol{\Xi}} \|\mathbf{B} - \mathbf{G}\boldsymbol{\Xi}\|_2 + \lambda\|\boldsymbol{\Xi}\|_0, \tag{13}$$

is solved by using a modified sequential thresholding least-squares algorithm [37], where a grid search is performed to find the optimal values for $\lambda$ in each component.

In WSINDy, the test function $\phi$ plays a key role in achieving its aim. To obtain $\phi$, the function used is based on the available code of WSINDy for ODEs [30] is of the form

$$\phi_{m,p}\left(\frac{t}{m\Delta t}\right) = \left(1 - \left(\frac{t}{m\Delta t}\right)^2\right)^p, \qquad (14)$$

where $m \in \mathbb{Z}^+$ is the support parameter such that $\phi$ is supported on $2m + 1$ points of spacing $\Delta t$, and $2p$ is the degree of $\phi$, where $p \in \mathbb{Z}^+$. The particular functional structure of $\phi$ has been selected for its smoothness, compact support, and ease of numerical evaluation. This combination of traits allows WSINDy to effectively average out high-frequency noise in the data. Additionally, these functions are separable or expressed as the multiplication of univariate functions, facilitating rapid computation of convolutions within the linear system. Importantly, they are obtained from the data itself, making them adaptable for application across a wide range of dynamical systems. In the following, a fast outline to construct the test functions is provided, see Eq. 4.3 and Appendix A from Ref. [37] for details:

a. Support and degree of $\phi$. The smoothness and decay properties of $\phi$ are influenced by $m$ and $p$, being these two factors crucial for the effectiveness of the method to guarantee decay in both real and Fourier spaces. In the real space, these two variables are related by

$$p = \min\left\{p: \phi\left(-1 + \frac{1}{m}\right) \leq \tau\right\},$$

then

$$\tau \leq \left(\frac{2m_d - 1}{m_d^2}\right)^p, \qquad (15)$$

where $\tau = 10^{-10}$ is a decay tolerance hyperparameter in real space.

b. Computation of wavenumbers. The data $\mathbf{Y}(t)$ is transformed into the frequency domain as the result of applying the fast Fourier transform to $\mathbf{Y}(t)$ $(\mathcal{F}(\mathbf{Y}(t)))$ along each component. Then, the wavenumbers $k^*$ are located at the corner where a $L_2$ fitting piecewise-linear approximation of two pieces of the cumulative sum of the data values in $|\mathcal{F}(\mathbf{Y}(t))|$ is the minimum.

c. Fitting $\phi$ with a Gaussian distribution. To achieve noise filtering, the test function $\phi$ is approximated by a Gaussian distribution, $\phi(t) \approx \rho_\sigma(t)$, with $\rho_\sigma(t) = \frac{1}{\sqrt{2\pi\sigma^2}} e^{-t^2/2\sigma^2}$, and $\sigma = \frac{m\Delta t}{\sqrt{2p+3}}$. Thus, to enforce the decay of $\hat{\phi}\left(\frac{2\pi}{N\Delta t} k_d^*\right)$, $m$ and $p$ are related by

$$\frac{2\pi}{N\Delta t} k_d^* = \hat{t}\hat{\sigma} = \frac{\hat{t}}{\sigma} = \hat{t}\frac{\sqrt{2p+3}}{m\Delta t}, \qquad (16)$$

where $\hat{t}$ is a rate of decay (user-defined) hyperparameter in the Fourier space to enforce $k^*$ falls $\hat{t}$ standard deviations into the tail of the spectra $\hat{\phi}$. Solving (13) and (14), the following expression is obtained

$$F(m) = \log\left(\frac{2m-1}{m^2}\right)\left(4\pi^2 k^{*2} m^2 - 3N^2 \hat{t}^2\right) - 2N^2 \hat{t}^2 \log(\tau). \qquad (17)$$

After numerically finding $m$, either (13) or (14) can be used to get $p$.

## 3. The proposed algorithm: WmSINDy

Despite its success, mSINDy method relies on the traditional SINDy method that computes the pointwise derivatives to the collected time-series data to initiate the algorithm, during the optimization, and in the sparse regression procedure. Calculating the derivatives of measurement data affected by noise presents a difficult task due to the randomness of the noise, which can result in high fluctuations and significant errors in the computation [38]. The mathematical formulation of numerical differentiation is typically ill-posed, making it difficult to choose the appropriate computational method and its parameters [27]. Even if the true signals are smooth, noise can make the data appear as if it has abrupt changes, leading to singularities or discontinuities characterized by spikes in the derivative or even worse. Furthermore, noise can lead to zones of the domain where the derivative does not exist, requiring more sophisticated methods [39].

These challenges necessitated the development of novel approaches to deal with noisy data such as the total variation regularized differentiation approach in the original SINDy [2], the Savitzky-Golay filter [25], least squares [40], splines [41], proximal gradient [42], Kalman smoothing [43], and others to improve the accuracy and reliability of numerical derivative calculations in the presence of noise. Although successful, smoothing noisy data may eliminate important features in the data and the effectiveness of smoothing often depends on the choice of parameters and the characteristics of the noise. Often, sophisticated smoothing techniques are computationally intensive, which makes them difficult to use in cases of large measurements. These severely limit the robustness of the algorithm to noise and hence its applicability. Here, we relax this assumption and use the weak formulation of the dynamics to circumvent this problem in mSINDy.

The weak formulation, as demonstrated in [30], allows the recovery of ODEs without relying on pointwise derivative approximations, which leads to better noise robustness. By setting the ODE into its weak form, the method avoids the approximation of pointwise derivatives, resulting in effective machine-precision recovery of model coefficients from noise-free data and robust identification of ODEs in the relatively large noise regime. This approach is particularly valuable as it enables the identification of a wide variety of ODEs under relatively large noise environments. In this regard, we propose to infuse the WSINDy method into the mSINDy method to avoid the calculation of the point-wise derivatives. The implementation of the weak formulation into the modified SINDy is possible by introducing Eq. 13 into Eq. 4 such that the $e_d$ is replaced by the residual error ($e_r$) as

$$\mathcal{L}_{WmSINDy}(\Xi, \widetilde{\mathbf{N}}) = e_r + e_s, \tag{18a}$$

$$\mathcal{L} = \mathcal{L}_{WmSINDy}(\Xi, \widetilde{\mathbf{N}}) = \left\|\mathbf{B}(\widetilde{\mathbf{X}}) - \mathbf{G}(\widetilde{\mathbf{X}})\Xi\right\|_2^2 + \sum_{j=q+1}^{N-q} \sum_{i=-q, i \neq 0}^{q} \omega_i \left\|\mathbf{y}_{j+i} - \widetilde{\mathbf{n}}_{j+i} - \widetilde{\mathbf{F}}^i(\widetilde{\mathbf{x}}_j)\right\|_2^2. \tag{18b}$$

In the proposed approach (depicted in Fig. 1), adapted to the scheme of [26], the estimated noise is subtracted from the measurement to obtain the denoised time series. Thus, in an optimization framework, we aim to simultaneously identify the parsimonious dynamical system $\mathbf{G}(\mathbf{X})\Xi$ responsible for generating the time series $\widetilde{\mathbf{X}}$, and noise $\mathbf{N}$. The optimization problem searches for the solution

$$\Xi, \widetilde{\mathbf{N}} = \arg\min_{\Xi, \widetilde{N}} \mathcal{L}(\Xi, \widetilde{\mathbf{N}}),$$

$$\text{s.t. } (|\Xi| < \lambda) = 0. \tag{19}$$

The WmSINDy framework is shown in Algorithm 1. Unlike mSINDy, in step 1, we first initialize the algorithm by calculating a $\Xi$ by WSINDy method. Then, in step 3, we replace the $e_d$ term in mSINDy by $e_r$ during the optimization. Finally, in step 8, during the update of $\Xi$, we perform the sparse regression such that we use the weak convolution form. In all cases, when simulating the systems, we use the **odeint** solver, with the relative and absolute tolerances set to $10^{-12}$, values that allow high numerical accuracy of the simulation of the considered dynamical systems, without compromising computational costs, and the time step $dt = 0.01$. To minimize the loss function $\mathcal{L}$, the Adam optimizer [44] with a learning rate of $0.001$ is used with a maximum of 5000 iterations for each loop[26]. All numerical experiments were developed in a Tensorflow environment [35] at the Holland Computing Center server [45]. It is worth noting that WSINDy computes the convolution between the test functions and the filtered data. In the WmSINDy, we avoid filtering the data to reduce the number of hyperparameters.

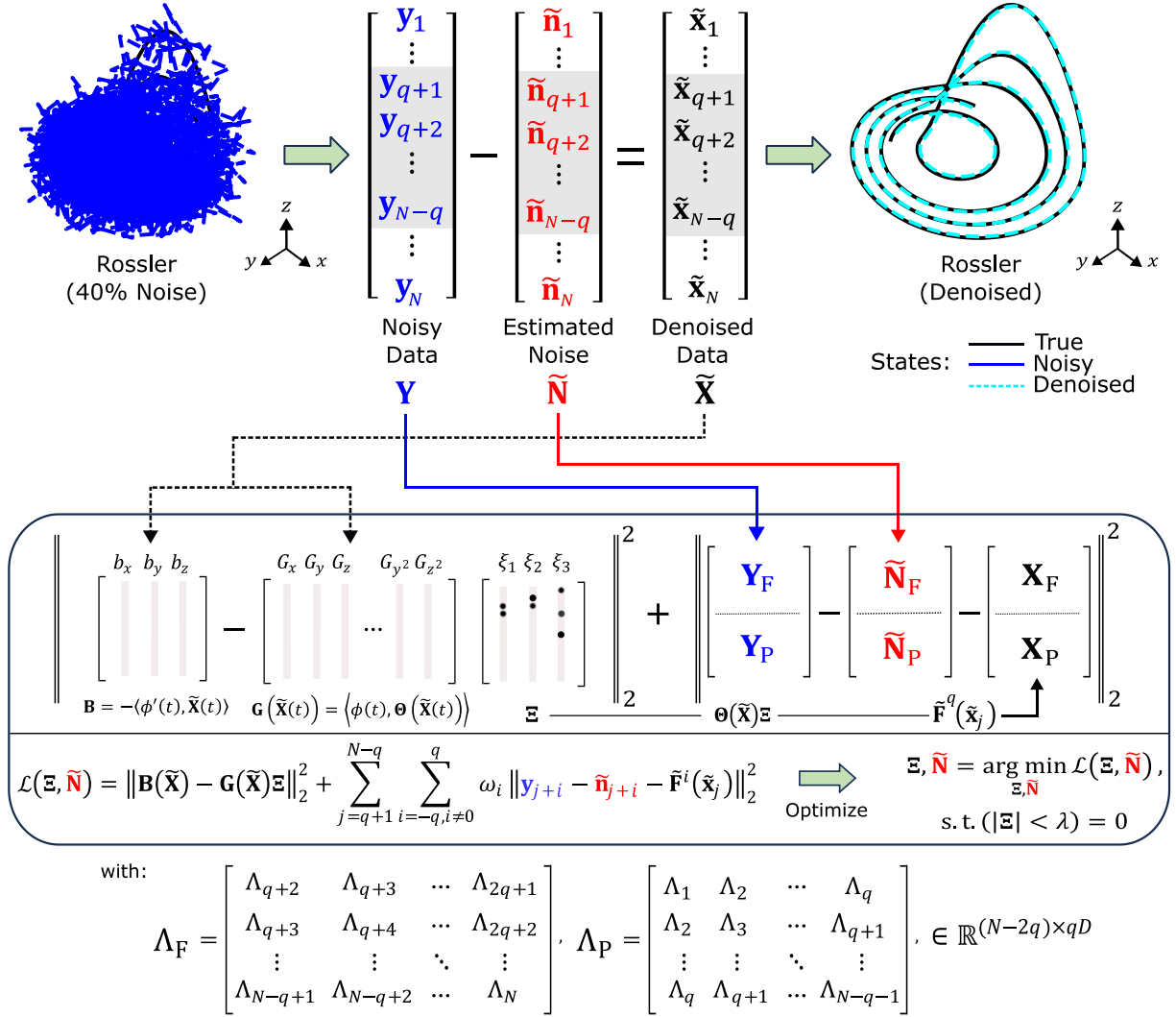

**Fig. 1.** The Weak modified SINDy algorithm

**Algorithm 1:** WmSINDy: Weak modified Sparse Identification of Non-linear Dynamics

**Input:** $Y(t), \Theta(*), dt, \lambda, N_{loop}, \omega, k, N, \tau, \hat{\tau}$.
**Output:** $\Xi, \widetilde{N}(t)$

1. $\Xi, \phi, \phi' = \text{WSINDy}(\Theta(Y(t)), k, \tau, \hat{\tau})$ in (8)  ▷ Initialize $\Xi$, and test functions $\phi$ and $\phi'$
2. $\widetilde{N}(t) = \text{zeros}(\text{size}(Y(t)))$  ▷ Initialize noise of the system

**for** $i = 1$ to $N_{loop}$ **do**  ▷ Denoising and learning system model
3.    Optimize     $\mathcal{L}(\Xi, \widetilde{N}(t))$ in (14).
4.    Apply sparsity     $(|\Xi| < \lambda) = \mathbf{0}$.
5.    Get new estimate of the state     $\widetilde{X}(t) = Y(t) - \widetilde{N}(t)$.
6.    From the new state $\widetilde{X}$, obtain     $\phi, \phi'$.
7.    Calculate     $B(\widetilde{X}(t))$ and $G(\widetilde{X}(t))$.
8.    Perform sparse regression     $(|\Xi| \neq 0) = G(\widetilde{X}(t)) \backslash B(\widetilde{X}(t))$
**end**

## 4. Performance comparison of the WmSINDy

### 4.1 Comparison with mSINDy and WSINDy methods

To reflect the accuracy of system identification and noise characterization, we study the robustness to noise, amount of data, thresholding parameter $\lambda$, prediction step $q$, and success rate (if all corresponding elements in $\Xi$ and $\widetilde{\Xi}$ are non-zeroes) [26]. The criteria to compare mSINDy and WmSINDy are presented in Table 1, while the normalized parameter error is used from the three methods: WSINDy, mSINDy, and WmSINDy.

**Table 1.** Criteria for comparison

| Error | Formula |
|---|---|
| Noise identification | $E_N = \dfrac{1}{N}\sum_{i=1}^{N}\|\mathbf{n}_i - \widetilde{\mathbf{n}}_i\|_2^2$ |
| Vector field | $E_f = \dfrac{\sum_{i=1}^{N}\|f(\mathbf{x}_i) - \tilde{f}(\mathbf{x}_i)\|_2^2}{\sum_{i=1}^{N}\|f(\mathbf{x}_i)\|_2^2}$ |
| Between forward simulation and true trajectory | $E_F = \dfrac{1}{\|X\|_F^2}\sum_{i=1}^{N-1}\|\mathbf{x}_i - \widetilde{F}^i(\mathbf{x}_1)\|_2^2$ |
| Normalized parameter | $E_p = \dfrac{\|\Xi - \widetilde{\Xi}\|_2}{\|\Xi\|_2}$ |

To assess the effectiveness of identification methods and conduct a comparison of the algorithms, the Lorenz attractor is considered:

$$\dot{x} = \sigma(y - x), \quad (20a)$$
$$\dot{y} = x(\rho - z) - y, \quad (20b)$$
$$\dot{z} = xy - \beta z, \quad (20c)$$

where $\sigma = 10$, $\rho = 28$, and $\beta = 8/3$. These equations describe a simplified system of atmospheric convection, which under these parameters the system exhibits chaotic behavior [46]. In WSINDy for ODEs [30], there are several hyperparameters to be tuned. In all cases, for both WSINDy and WmSINDy, we keep the same values based on the existing codes for the WSINDy algorithm, except for the number of query points, $\hat{\tau}$, and the scale factor. The number of query points is set equal to the number of sampling points, $\hat{\tau} = -2$ [32] and the scale factor to zero because there is no considerable difference between the biggest and smallest values. In the following, we perform a comparison between the baselines and the proposed method by considering the effects of noise level, data length ($N$), threshold ($\lambda$), prediction step ($q$), and iteration loops ($N_{loops}$).

### 4.1.1 Effect of noise

Various levels of Gaussian noise, from 0 to 50%, are introduced to create the noisy training data. The noise level is expressed as [26]

$$\text{Noise Level}(\%) = \sqrt{\frac{\text{var(Noise)}}{\text{var(Signal)}}} \times 100\% = \frac{\text{std(Noise)}}{\text{std(Signal)}} \times 100\%. \tag{21}$$

At every noise level, 10 distinct sets of noisy data are created and employed as input for the three methods. For both mSINDy and WmSINDy methods the hyperparameters, $N_{loop} = 6$, $\lambda = 0.2$, and $q = 1$ are used. Unless otherwise stated, the library is built using polynomial terms up to second order, excluding the constant term. The simulation of the Lorenz attractor starts from the initial condition [5, 5, 25], and runs for a time of 25 seconds. For the simulated trajectory, the identified model is computed for 6 seconds forward in time[26]. In the figures where comparisons are made, the black circle denotes the median from ten runs, while the violin shape illustrates the spread of errors. The height of each shape represents the range of the result and the width at any point represents the density of data points at that value.

In Fig. 2(a), for noise identification error, both mSINDy and WmSINDy show an increasing error as the noise intensity increases, nonetheless, except for noise-free cases WmSINDy maintains lower error compared with mSINDy. Note that mSINDy displays a steep increase in error as noise levels increase, with a significant jump around 5-10% noise level. It is worth mentioning that in the available code of mSINDy, the $\lambda$ parameter can be adaptively tuned depending on the noise level; however, this search may be difficult to perform since the noise intensity may not be known a priori. For vector field error (Fig. 2(b)), except for the noise-free case, WmSINDy performs better than mSINDy by showing lower errors and less variability, and also WmSINDy performs better than WSINDy by showing lower errors for the medians in the high noise intensity region. For short-term prediction error (Fig. 2(c)), in general, WmSINDy exhibits very low errors, especially in the case of 5-10% noise level, and it keeps lower errors than mSINDy through different noise level intensities. For the normalized parameter error (Fig. 2(d)), in the noise-free case, the WSINDy

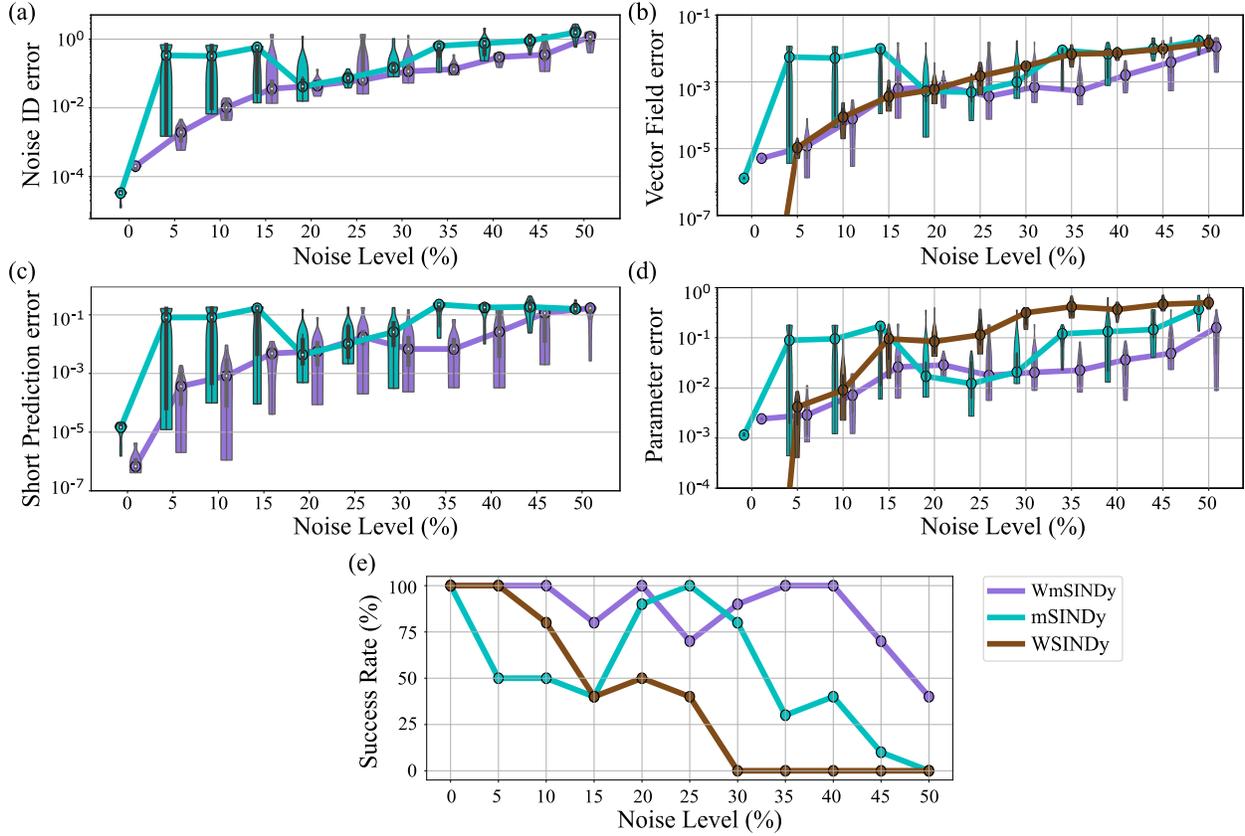

**Fig. 2.** Comparison for different noise levels.

(error: $4.54 \times 10^{-13}$) considerably outperforms both the mSINDy and the WmSINDy, and it is possible in these cases that WSINDy can recover the $\Xi$ coefficients with a precision that is only limited by the capabilities of the computer used [30]. Except for the noise-free case, WmSINDy demonstrates lower parameter errors and remains relatively stable compared to mSINDy and WSINDy, while mSINDy shows higher error and more variability compared to WmSINDy, and from 20% onwards, WmSINDy exhibits the highest errors. For the success rate (Fig. 2(e)), WmSINDy has a higher success rate across different noise levels, indicating more consistent and reliable performance. The success rate drops slightly with increasing noise but remains relatively high up to 45%. mSINDy displays a lower and more variable success rate compared to WmSINDy. Finally, WSINDy shows the lowest success rate, especially at higher noise levels. For the three methods, the success rate drops significantly at higher noise levels, nonetheless, the WmSINDy stands out over the other two. Note that in Figs. 2(a) and (c), only methods mSINDy and WmSINDy are compared, as WSINDy is not applicable in these specific cases. Specifically, the WSINDy's process is inherently unable to characterize noise and simulate the states, and therefore cannot be included in these analyses. As such, we will not include WSINDy in any subsequent plots that compare the methods for the noise identification. These results show the interesting combination of both mSINDy and WSINDy, in which, especially, in the normalized parameter error, the improvement to both methods is more visible in the high noise regime. Among the methods evaluated, WmSINDy emerges as the best approach. It outperforms the other methods in terms of reliability and effectiveness. mSINDy performs adequately at noise-free but struggles in low noise, if $\lambda$ is not properly set, and in high noise, showing significant increases in errors and decreases in success rate. For this numerical experiment, when using WmSINDy, we found the best results for $\lambda$ in the range of [0.2, 0.4]. A value below it produces overfitting, while setting a higher value we obtain underfitted models.

### 4.1.2 Effect of data length

We tested the effectiveness of the three methods, WSINDy, mSINDy, and WmSINDy with different amounts of available data, keeping the noise level (40%) the same. For the 10 distinct sets of noisy, the hyperparameters, $N_{loop} = 6$, $q = 3$, and $\lambda = 0.2$ are used, with $x_0 = [-5, 5, 25]$. For the short-term prediction error, 24% of the data points are used, which corresponds to percentage used in the previous section.

For noise identification errors (Fig. 3(a)), WmSINDy tends to have lower errors compared with mSINDy. Also, the error in WmSINDy decreases as the data length increases, while for mSINDy the error seems to have a non-decreasing behavior. For vector field error (Fig. 3(b)), WmSINDy maintains lower error consistently across different data lengths compared with mSINDy and WSINDy. For short-term prediction errors (Fig. 3(c)), WmSINDy exhibits very low errors, especially for shorter data lengths; mSINDy shows low errors as well, but higher errors compared with the proposed method. Note that WSINDy is excluded from Figs. 3(a) and (c) because it cannot identify the noise. For parameter error (Fig. 3(d)), from the three methods, WmSINDy demonstrates the lowest errors, and it tends to decrease as the data length increases, indicating a good parameter estimation accuracy, while mSINDy performs similarly to WmSINDy but with higher errors. Finally, WSINDy shows consistently high errors and less improvement with the data length increase. For the success rate (Fig. 3(e)), the WmSINDy reaches 90% in the system identification, at 1500, 2250, and 2500 data points. In contrast, the maximum value achieved by mSINDy is 40% when 1250 points are used. Thus, WmSINDy has a higher success rate through different data lengths, which indicates that it is more consistent and reliable. In summary, WmSINDy stands out as the most successful option, while mSINDy achieves satisfactory results, it falls short of WmSINDy's performance. WSINDy comes in last place in terms of effectiveness.

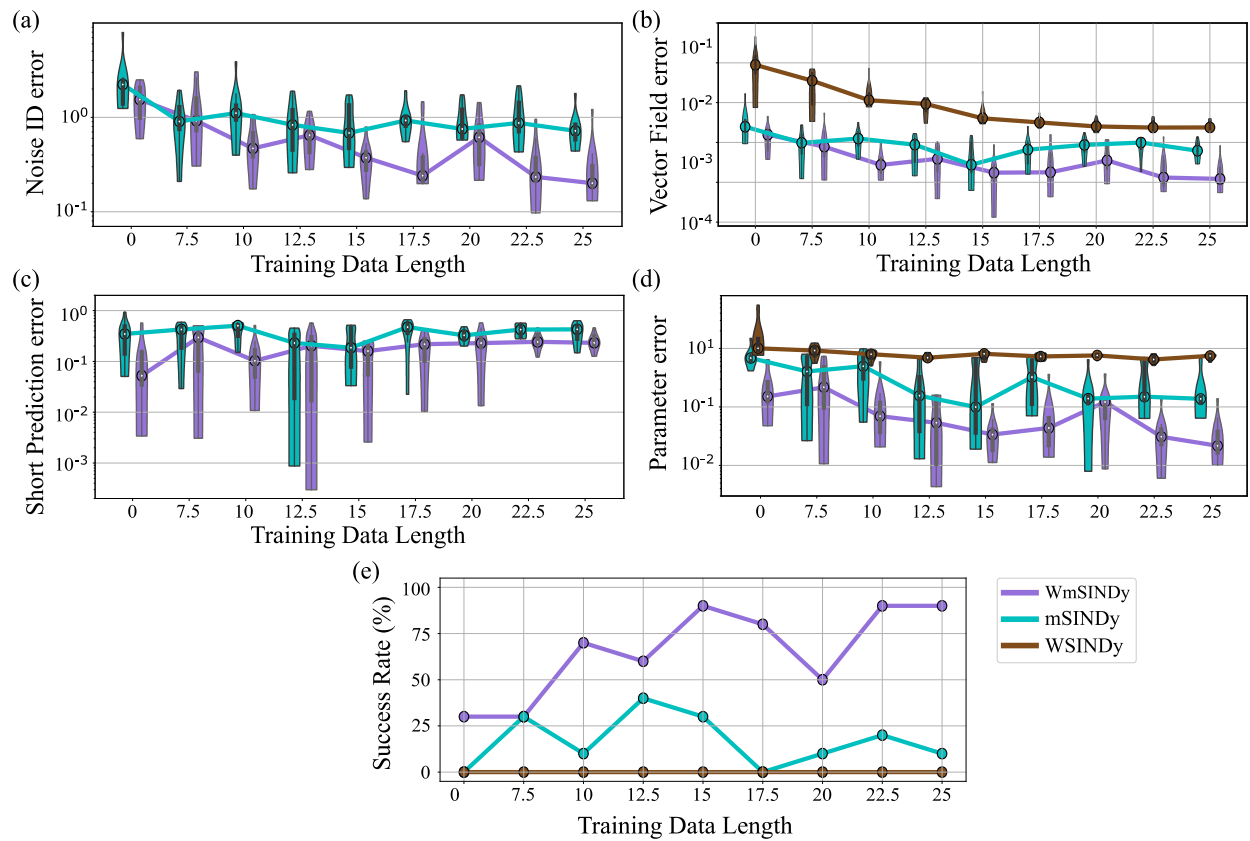

**Fig. 3.** Comparison for different training length data

### 4.1.3 Effect of training thresholding parameter $\lambda$

When $\lambda$ is small, the model may become overfitted and include many unnecessary terms, which may produce a stagnancy in a local minimum. While, if the $\lambda$ value is large, the model may become underfitted and exclude many relevant terms. For this numerical experiment, the Lorenz attractor is simulated with $x_0 = [-5, 5, 25]$ for a time unit of 25 with 40% of noise. At every $\lambda$ value, 10 distinct sets of noisy data are used, and the numerical experiment is run for, $N_{loop} = 6$, and $q = 1$ for both mSINDy and WmSINDy methods.

For noise identification error (Fig. 4(a)), WmSINDy generally maintains lower errors across different $\lambda$ values and less variability compared to mSINDy. This indicates that when $\lambda$ is not accurately known, WmSINDy is slightly more robust in noise identification. For parameter error (Fig. 4(b)), the optimal region of $\lambda$ values is from 0.2 to 0.5 for both methods. mSINDy displays slightly lower errors than WmSINDy but has much higher degree of variability. For vector field error (Fig. 4(c)), WmSINDy performs better than mSINDy, particularly at $\lambda$ values around 0.3 to 0.6. This indicates the WmSINDY is more accurate at estimating the vector field compared to mSINDy. The comparison for the success rate is shown in Fig. 4(d), which shows that, for $\lambda \in [0.05, 0.1]$ and $\lambda \in [0.8, 1]$, both methods are not able to accurately identify the dynamics correctly. For mSINDy, the best case is achieved when $\lambda = 0.2$ resulting in a success rate of 80%, whereas the best cases for WmSINDy are 100% success rates for $\lambda = 0.3$ and $\lambda = 0.4$. In this sense, WmSINDy will outperform mSINDy for a high noise regime if the $\lambda$ hyperparameter is correctly chosen. In summary, WmSINDy appears to be a more reliable and robust method for different $\lambda$ values, particularly, at $\lambda = 0.3$ and 0.4 where the success rate is 100%, while mSINDy struggles with higher errors and higher variabilities.

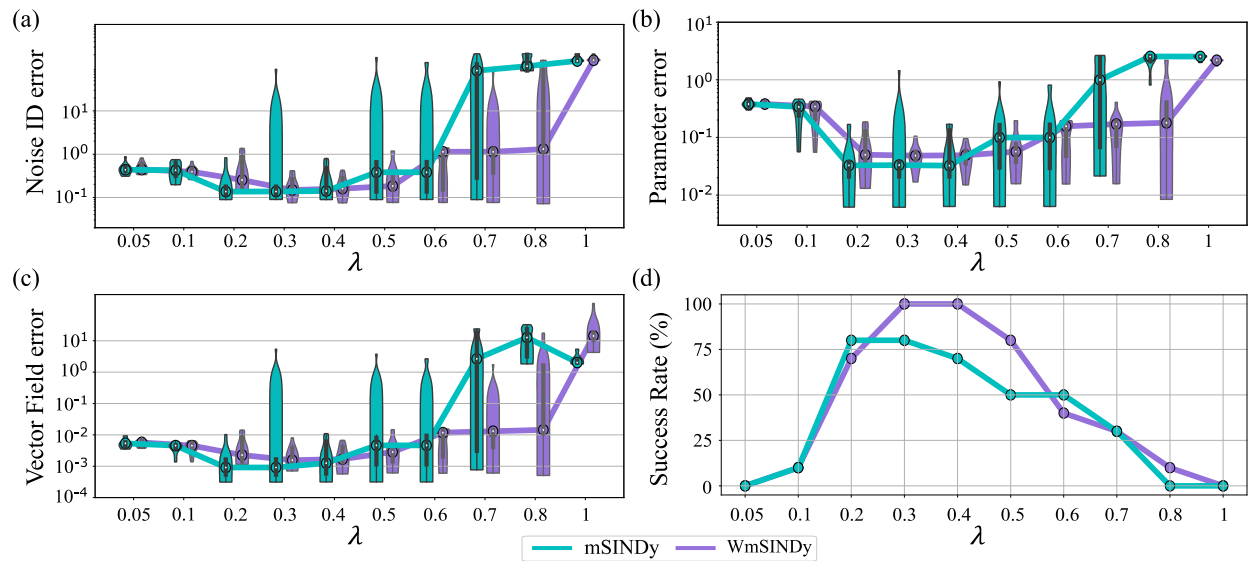

**Fig. 4.** Comparison for different $\lambda$ values

### 4.1.4 Effect of prediction step $q$

The prediction step $q$ controls how many timesteps forward and backward the identified model is simulated for when computing the simulation error term, $e_s$. A larger $q$ leads to simulating the identified model over a longer time horizon. To study the effect of $q$ on the performance of mSINDy and WmSINDy, we use a time of 25 s, $x_0 =[-5,5,25]$, $N_{loop} = 6$, $\lambda = 0.2$, and noise at 40%, and depict the results of the study in Fig. 5. Note that we do not compare with WSINDy in this subsection $q$ is not a parameter in that method.

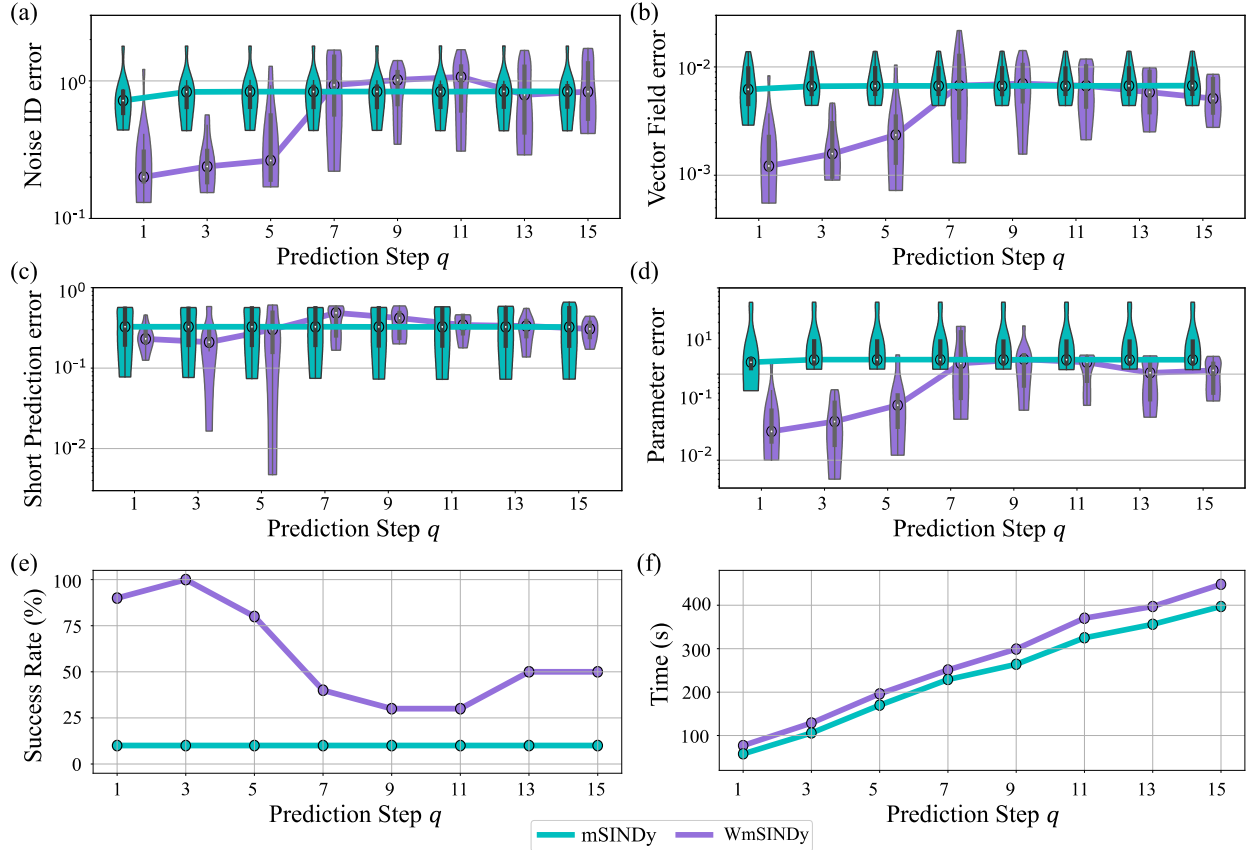

**Fig. 5.** Comparison for different prediction step $q$.

The results for noise identification error are presented in Fig. 5(a), which shows that, when $q$ is small, the noise characterization made by WmSINDy is the best. On the other hand, mSINDy is only minimally affected by changes in $q$. Furthermore, WmSINDy generally results in lower error through different prediction steps $q$ compared to mSINDy. When $q = 7$, the error produced by WmSINDy increases to match the behavior of the mSINDy, indicating that WmSINDy performs best with small values of $q$. For the vector field error, shown in Fig. 5(b), WmSINDy produces consistently smaller errors for different prediction steps, while for mSINDy the results have little dependence on the prediction step $q$ except for $q = 1$. For $q \geq 7$, WmSINDy produces similar results to mSINDy, but with slightly lower errors. For the short-term prediction error (Fig. 5(c)), mSINDy exhibits a stable error that is lower than WmSINDy, except from $q = 7$ onwards. Nonetheless, the WmSINDy provides better performance in the range 1-3, which is also desirable for the computational time standpoint. For the parameter error (Fig. 5(d)), WmSINDy results in lower errors than mSINDy for all $q$ values. For $q \geq 8$, WmSINDy produces a constant median just like mSINDy, but with lower errors. The success rate is presented in Fig. 5(e), which shows that WmSINDy has a significantly higher success rate, especially for small $q$ values, with a peak success rate of 100% at $q = 3$. The success rate decreases slightly at higher prediction steps but remains higher than mSINDy. Regardless of the value of $q$, mSINDy consistently results in a low success rate at around 10%, indicating that $q$ has minimal effect on the success of mSINDy and that mSINDy has limited overall performance. We present the calculation time of one experiment in Fig. 5(f), which demonstrates that the required computation time increases for both methods as $q$ increases. However, mSINDy is less computationally expensive than WmSINDy regardless of the value of $q$. WmSINDy requires more calculation time because of the convolution operations to get **G** and **B**, and this increased time is minimal compared to mSINDy. In summary, WmSINDy generally performs better across all metrics, resulting in lower errors and higher success rates. It is preferable to use a small $q$, since it is going to strongly influence the running process. In the following, in most of the cases we are going to use $q = 1$.

### 4.1.5 Effect of optimization iteration $N_{loop}$

The $N_{loop}$ value controls the number of iterations performed by the algorithm during the identification. If this parameter is too small, then the algorithm will be unable to reach an optimal minimum and will fail to identify the dynamics or characterize the noise properly. Large values of $N_{loop}$ provide better accuracy for the price of considerable increases in the computational cost. To study this parameter, we simulate the Lorenz attractor with an initial condition of [5,5,25] for a time-unit of 25 with a 40% of noise. We consider $N_{loop} \in \mathbb{Z} : N_{loop} \in [1,8]$ and apply mSINDy and WmSINDy to ten distinct sets of noisy data for $\lambda = 0.2$ and $q = 1$. In this case, we exclude WSINDy because the $N_{loop}$ hyperparameter is not part of that method.

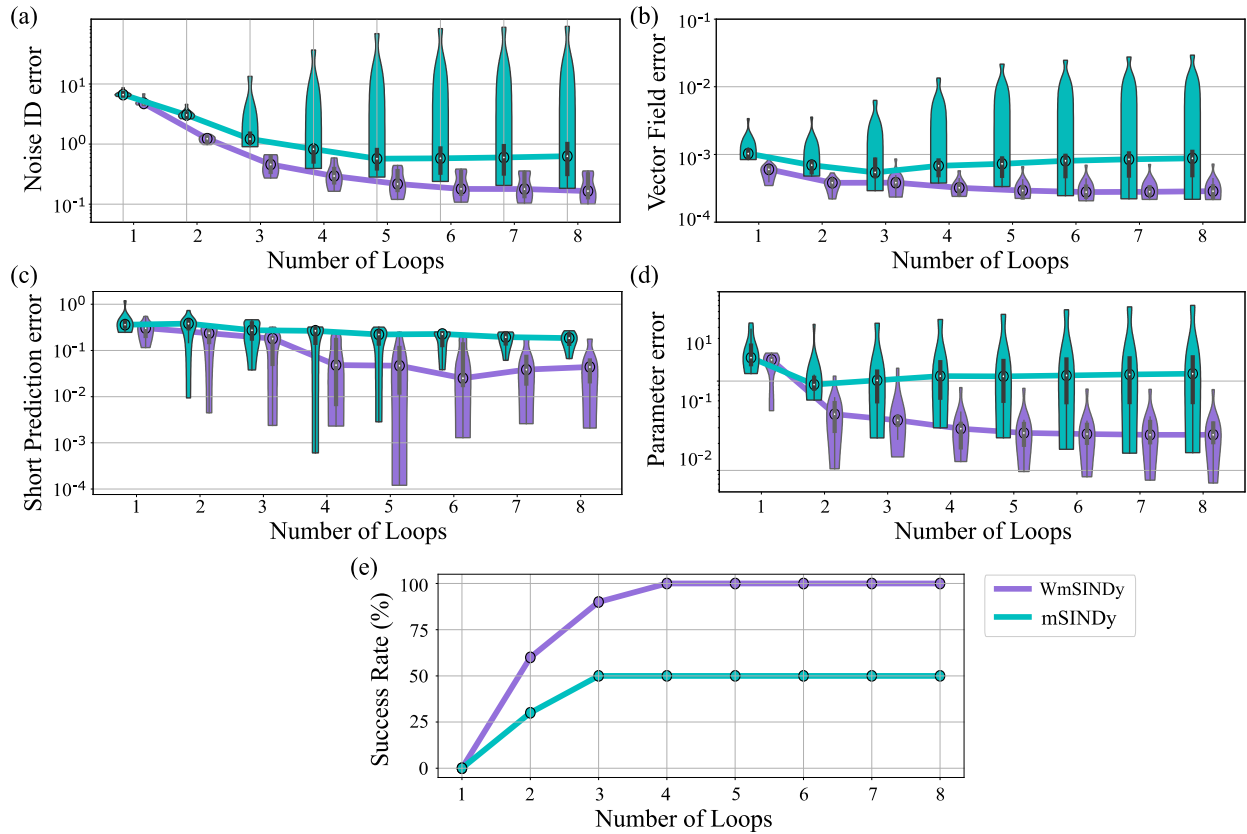

**Fig. 6.** Comparison for different number of loop values

For WmSINDy, the errors shown in Figs. 6(a)-(d) decrease as the $N_{loop}$ increases, and the results tend to converge for $N_{loop} \geq 5$. A similar behavior is observed for mSINDy; however, the ranges in Figs. 6(a,b,d), the mean and the maximum errors are considerably high. As a result, the reliability of mSINDy is not guaranteed depending on the noise. Note that in Figs. 6(c), WmSINDy results in higher error ranges, but shows great performance at the lower bounds. In all cases, WmSINDy provides better results than mSINDy. Figure 6(e) displays that WmSINDy has consistently higher success rates compared to mSINDy. When $N_{loop} = 4$, WmSINDy converges to a success rate of 100%, while mSINDy converges to 50% success rate at this value. Overall, while both methods benefit from a greater value of $N_{loop}$, the results indicate that WmSINDy consistently outperforms mSINDy in terms of lower errors, less variability, and higher success rates for different numbers of loops evaluated in this analysis. Based on these results, we will use the $N_{loop}$ at least a value of 5 in the following, which is a reasonable value to characterize the noise and identify the dynamics of the considered systems at a moderate time.

## 4.2 Comparison with mSINDy

The performance of the proposed algorithm is also demonstrated on a variety of canonical systems that exhibit regular and chaotic responses as well as limit cycle oscillations. Table 2 shows the parameters of the system and the hyperparameters for the optimization to be used. By studying these systems, we demonstrate the versatility of the proposed WmSINDy across a wide range of response types.

**Table 2.** Parameters of the system and the hyperparameters used in the canonical systems.

| System | Noise intensity | Initial conditions | Run time | $N_{loop}$ | $\lambda$ | $q$ |
|---|---|---|---|---|---|---|
| Rössler attractor | 40 % | [3,5,0] | 25 s | 6 | 0.08 | 3 |
| Lorenz 96 | 40 % | [1,8,8,8,8,8] | 25 s | 8 | 0.2 | 1 |
| Van der Pol oscillator | 40 % | [−2,1] | 10 s | 8 | 0.15 | 1 |
| Duffing oscillator | 40 % | [−2,2] | 25 s | 5 | 0.05 | 1 |
| Cubic oscillator | 20 % | [0,2] | 25 s | 5 | 0.08 | 1 |
| Lotka-Volterra | 30 % | [1,2] | 10 s | 5 | 0.2 | 1 |

### 4.2.1 Rössler attractor

The equation that governs this system is

$$\dot{x} = -y - z, \quad (22a)$$
$$\dot{y} = x + ay, \quad (22b)$$
$$\dot{z} = b + z(x - c), \quad (22c)$$

where $a = 0.2$, $b = 0.2$, and $c = 5.7$. This system was designed for studying chaos, intended to be easier to analyze than the Lorenz attractor [47]. The performance of the two algorithms is shown in Fig. 7. WmSINDy has significantly lower errors compared to mSINDy across the different error categories shown on the x-axis. Specifically, the differences in the medians of the errors range from two to four orders of magnitude. The bar chart on the right shows the success rates, where WmSINDy achieves a success rate of 60%, while method mSINDy has a much lower success rate of 30%. Overall, for this system, WmSINDy is more accurate and successful than mSINDy.

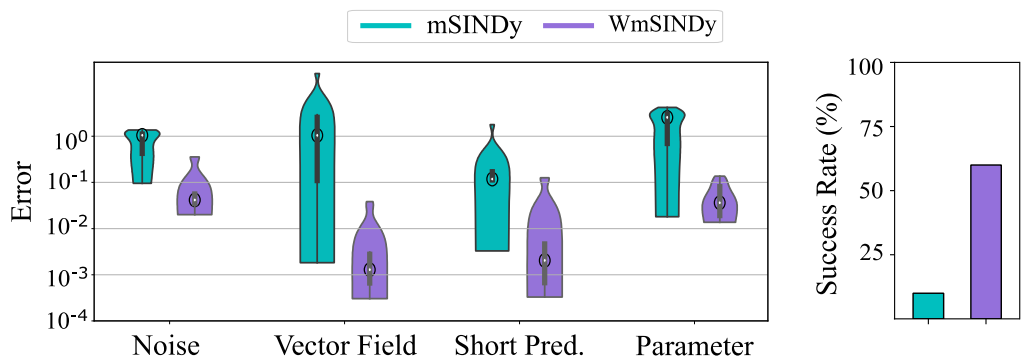

**Fig. 7.** Performance of the algorithms for the Rossler attractor.

### 4.2.2 Lorenz 96 model

The equation is given by:

$$\dot{x}_i = (x_{i+1} - x_{i-2})x_{i-1} - x_i + F, \tag{23}$$

for $i = 1,2,...,S$, with $S = 6$, and by defining $x_{-1} = x_{S-1}$, $x_0 = x_S$, $x_1 = x_{S+1}$, the model has 6 states. This model was developed for weather forecasting [48]. The system behaves chaotically when the forcing term is set to 8 and we simulate its response for a time of 25. For the short-term prediction, we use a time of 6. A library containing candidate functions is formed using polynomial terms up to the second order (including the constant term), resulting in 28 candidates. Fig. 8 shows that mSINDy excels in short-term prediction accuracy. Whereas WmSINDy performs better, with up to more than one order of magnitude of difference, in terms of noise identification, vector field, and parameter estimation errors, and has a significantly higher success rate, suggesting higher overall reliability when using WmSINDy, for this system. Note that in [26], 84 functions were used; however, despite using fewer functions here, the results are better. Overall, for the three studied systems that exhibit chaotic behavior, WmSINDy demonstrates significantly lower error rates and higher success rates than mSINDy. Thus, for this complex systems, the proposed WmSINDy captures better the underlying dynamics in high noise regimes.

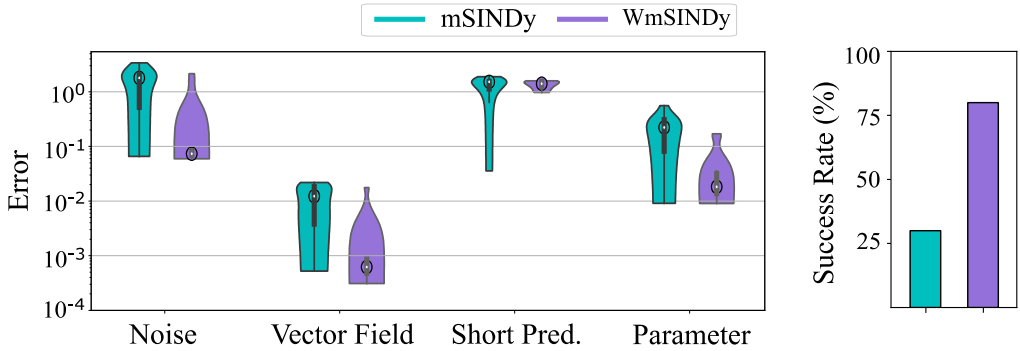

**Fig. 8.** Performance of the algorithms for the Lorenz 96 attractor.

### 4.2.3 Van der Pol oscillator

The equations that govern this system are

$$\dot{x} = y, \tag{24a}$$
$$\dot{y} = \mu(1 - x^2)y - x, \tag{24b}$$

where $\mu = 0.5$. These equations represent several important physical phenomena across various disciplines, including electrical systems to capture the behavior of relaxation oscillations[49], biological systems to model heartbeat rhythms and cardiac electrical activity[50], mechanical systems, to represent nonlinear vibrations and exited oscillations[51]. From Fig. 9, in noise, vector field, and prediction errors, mSINDy exhibits lower noise errors. While for parameter error, WmSINDy shows a lower error, suggesting better accuracy in parameter estimation. In the case of success rate, both methods achieve the same result. In a nutshell, for this system, mSINDy displays a more reliable choice for most of the metrics.

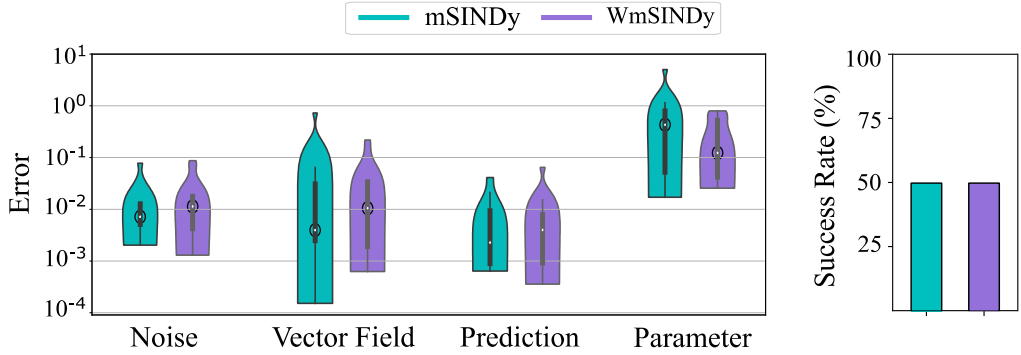

**Fig. 9.** Performance of the algorithms for the Van der Pol oscillator.

### 4.2.4 Duffing oscillator
The equations that govern this system are

$$\dot{x} = y, \tag{25a}$$
$$\dot{y} = -p_1 y - p_2 x - p_3 x^3, \tag{25b}$$

where $p_1 = 0.2$ and $p_2 = 0.1$, and $p_3 = 1$. This system represents the oscillations of a mass attached to a nonlinear spring and a linear damper[52,53]. Also, this model is effectively utilized in secure communications systems[54]. The comparison of Fig. 10 suggests that WmSINDy outperforms mSINDy across all metrics, especially because of the clear higher success rate and lower short-term predictions error.

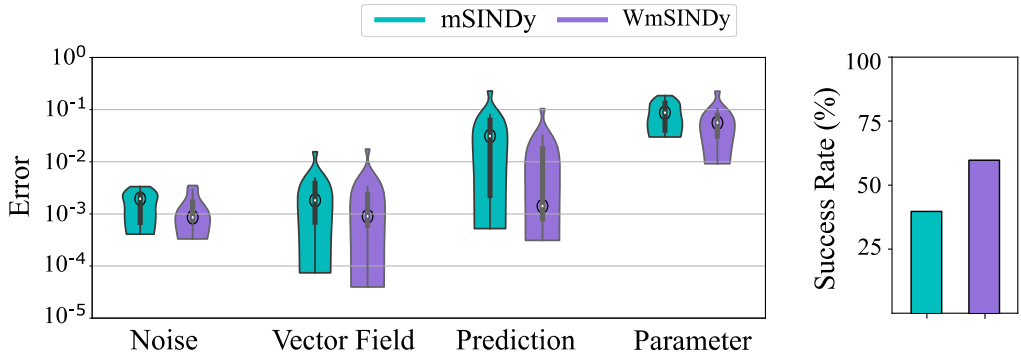

**Fig. 10.** Performance of the algorithms for the Duffing Oscillator.

### 4.2.5 Cubic oscillator
The equations that govern this system are

$$\dot{x} = p_1 x^3 + p_2 y^3, \tag{26a}$$
$$\dot{y} = p_2 x^3 + p_4 y^3, \tag{26b}$$

where $p_1 = -0.1$ and $p_2 = 2$, $p_3 = -2$, and $p_4 = -1$. This system is typically found in quantum mechanics to describe a particle within an anharmonic potential well[55]. In Fig. 11, overall, the medians of the errors are similar. However, some results provided by WmSINDy are considerably better than the ones provided by mSINDy, e.g., for the short-term prediction error we have up to two orders of magnitude of difference. In the case of success rate, both methods provide almost the same result, about 30%. For this system, the proposed WmSINDy could be a better option to perform system identification.

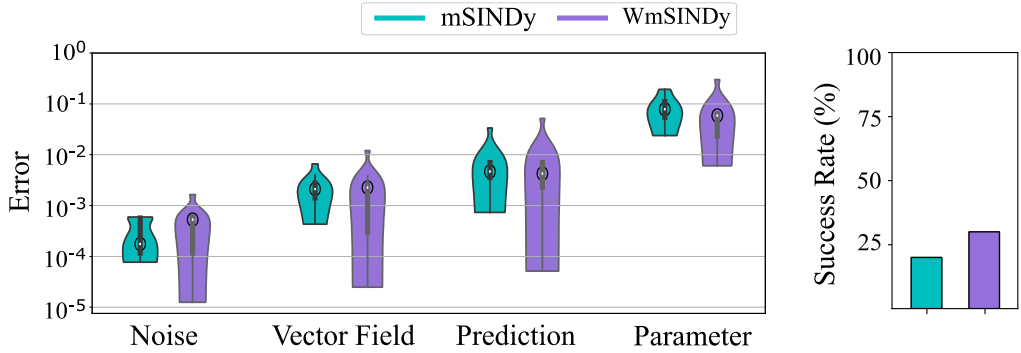

**Fig. 11.** Performance of the algorithms for the Cubic Oscillator.

### 4.2.6 Lotka-Volterra

$$\dot{x} = p_1 x - p_2 xy, \qquad (27a)$$
$$\dot{y} = p_2 xy - 2p_1 x, \qquad (27b)$$

where $p_1 = 1$ and $p_2 = 0.5$. These equations describe the dynamics of biological systems in which two species interact, typically as predator–prey[56]. In Fig. 12, for each metric, the medians of the errors are similar. The lowest values of the vector field error provided by mSINDy are considerably lower than the WmSINDy, however, the variance for mSINDy is notably bigger than mSINDy, and the same big variances for the other error metrics. Interestingly, the WmSINDy provides a 100% success rate compared with 60% of the mSINDy. Since the proposed method provides less variance in errors and a higher success rate, it could be a good option to deal with this system.

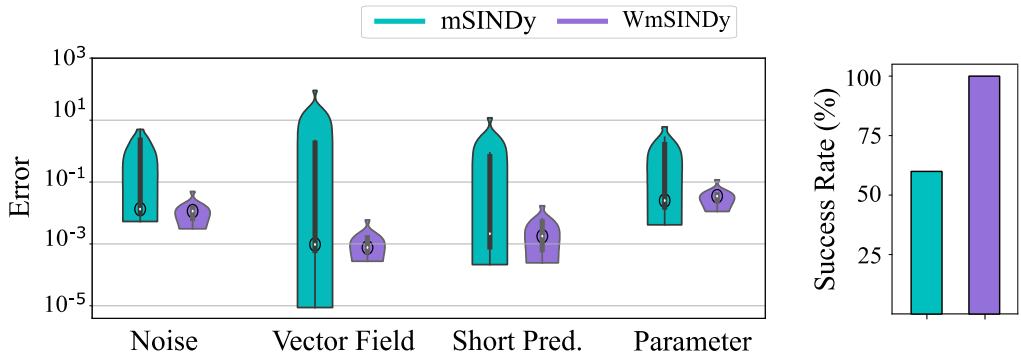

**Fig. 12.** Performance of the algorithms for the Lotka Volterra.

## 4.3 Discrepancy modeling

In scenarios where some aspects of the dynamics are known, it is convenient to introduce this knowledge to guide the system identification. Assume the observed dynamics in time-series data is expressed as

$$\dot{\mathbf{x}}(t) = \boldsymbol{f}(\mathbf{x}(t)) = \boldsymbol{g}(\mathbf{x}) + \boldsymbol{\Theta}(\mathbf{x})\Xi, \tag{28}$$

where $\boldsymbol{g}(\mathbf{x})$ and $\boldsymbol{\Theta}(\mathbf{x})\Xi$ are the known and missing dynamics, respectively. Discrepancy modeling is an approach that addresses the mismatch between the derivative of a known model and the observed dynamics of time-series data [57]. By combining mSINDy and discrepancy modeling, it was possible to identify the missing dynamics that were not captured by the known model.

To investigate the performance of the proposed algorithm, we consider the following system

$$\dot{x} = -10x + 10y + xy, \tag{29a}$$
$$\dot{y} = 28x - xz - y + 3z, \tag{29b}$$
$$\dot{z} = xy - 8/3\, z. \tag{29c}$$

where $\boldsymbol{g}(\mathbf{x})$ is given by

$$\dot{x} = -9.5x + 10.5y, \tag{30a}$$
$$\dot{y} = 27.6x - 1.1xz - 0.9y, \tag{30b}$$
$$\dot{z} = 1.05xy - 2.6z. \tag{30c}$$

For our comparison, the simulated system (Eq. 29) starts with the initial condition [5,5,25], runs for a time of 30, and noise with 40% intensity is added to the response. The hyperparameters are set as: $N_{loop} = 5$, $\lambda = 0.4$, $q = 4$. The identified model is simulated for a time of 7. The candidate functions library contains 9 candidate functions made up of polynomial terms up to the second order. Figure 13(a) shows that WmSINDy consistently outperforms mSINDy across all tasks in this comparison. WmSINDy shows lower errors for all four specific tasks and a significantly higher success rate. This suggests that when knowing some information of the system WmSINDy is generally more effective to characterize this complex system.

Although the goal of discrepancy modeling is to improve the accuracy and robustness of an existing model by incorporating the missing physics terms, it could be interesting to have an approach that can deal directly with a complex system, as presented in Eq. 29, without any prior knowledge. In this regard, without using an existing model, we apply both the mSINDy and the proposed WmSINDy. Figure 13(b) shows a comparison between mSINDy and WmSINDy. The proposed method excels in all metrics. Most importantly, it is highly reliable with a 100% success rate, making it the preferable choice for practical applications when no prior information is known. For this complex system, when there is no prior knowledge or when do exist known information, the WmSINDy shows better robustness than mSINDy in characterizing this system.

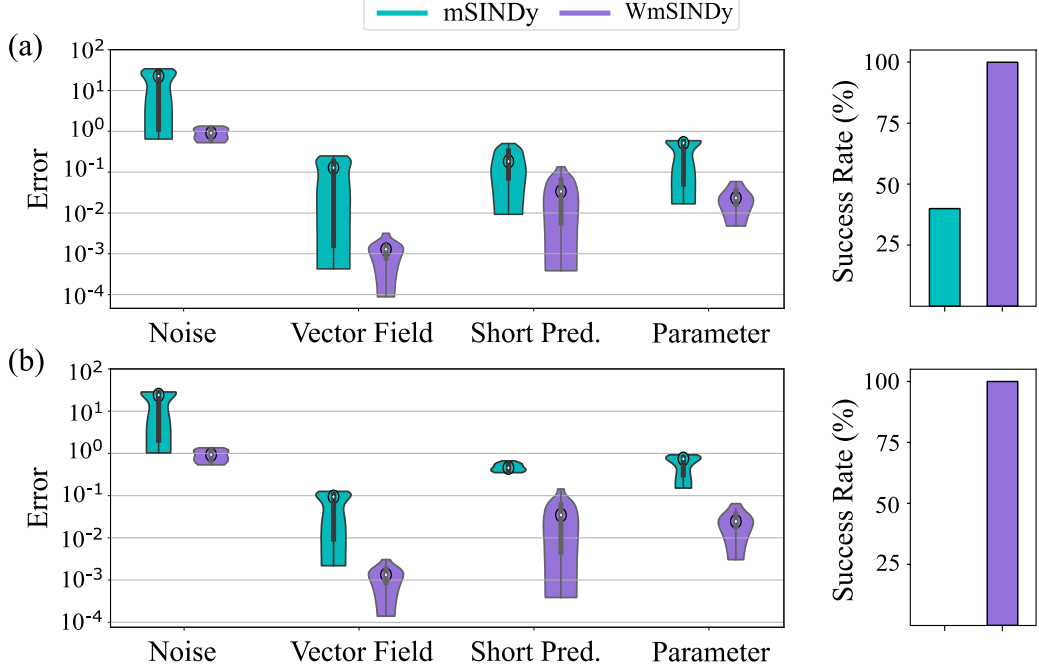

**Fig. 13.** Comparison between mSINDy and WmSINDy. (a) using model discrepancy; (b) without using modeling discrepancy.

### 4.4 Removing $e_r$ term

Kaheman et al. found in [26] that dropping the derivative error term from the cost function (Eq. (4)) does not significantly change the result. Instead, they found that dropping the derivative error term improves all the metrics. This improvement is believed to occur because dropping the derivative error term avoids the numerical derivative approximation error. Similarly, we drop the $e_r$ to know the impact of the WmSINDy. The numerical experiment is carried out by simulating the Lorenz attractor with and initial condition of [5,5,25], for a time of 25, under 40% of noise intensity. For both mSINDy and WmSINDy, the hyperparameters $N_{loop} = 8$, $q = 1$, and the library consists of second-order polynomials.

We present a comparison of the two methods. Figure 14 shows the results by using Eqs. (5) and (18) for mSINDy and WmSINDy, respectively. As can be seen, WmSINDy (dark purple) performs better than mSINDy (dark cyan) in all the metrics shown. By removing the $e_d$ term in mSINDy (light cyan), the performance of the method slightly improves compared to mSINDy when using the complete Eq. (5). Contrarily, removing $e_r$ (light purple) decreases the performance of WmSINDy compared to WmSINDy when using the complete Eq. (18), as can be seen that the errors have increased and the success rate decreased. For this system, it is better to omit the $e_d$ term in mSINDy to improve performance. Whereas for WmSINDy, the $e_r$ term that depends on the weak form of the dynamics, do not rely on pointwise approximations of the derivatives of the data. Thus, the results including the $e_r$ term show better performance because the weak form allows a more efficient handling of noisy data.

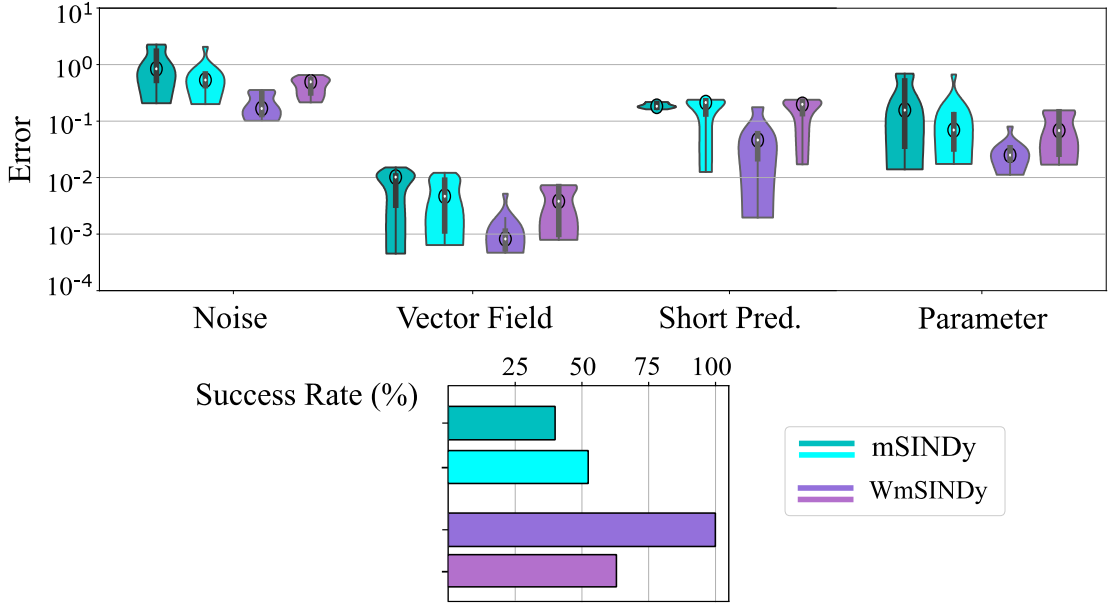

**Fig. 14.** Dark colors using Eqs. (4) and (18), (b) Light colors dropping $e_d$, $e_r$ term from Eqs. (4) and (18), in mSINDy and WmSINDy, respectively.

### 4.5 Learning non-zero mean noise

The numerical experiments so far have been focused on dealing with zero-mean noise distribution, such as Gaussian noise. As shown in [26], handling non-zero mean noise distribution is more challenging than handling Gaussian noise because non-zero mean noise introduces bias into the data. In this scenario, [26] proposed an iterative learning approach, which adapting to our algorithm consist in the following steps:

1. Learn the noise distribution by applying WmSINDy.
2. Update the data by subtracting the mean of learned noise.
3. Repeat step 2 to characterize the noise and identify the model.

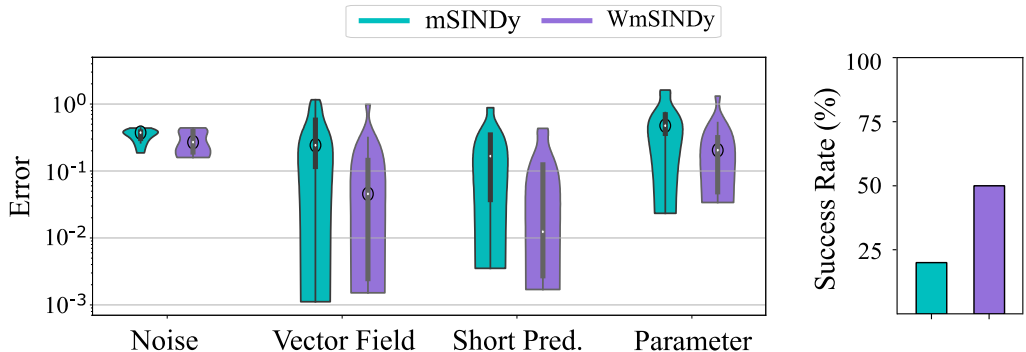

**Fig. 15.** Comparison between mSINDy and WmSINDy when learning non-zero mean noise.

For this study, step 2 is applied 3 times as done in the available code of [26]. The Van der Pol oscillator, serving as a standardized benchmark, is used to compare the performance of both mSINDy and WmSINDy methods. The system is the same as in Sec. 4.2.3, except that gamma noise with intensity of 30% is used. The gamma noise is a type of statistical noise that follows the gamma distribution. It is always positive and typically left-skewed, which results in a non-zero mean distribution. It is generally generated in fission events in nuclear reactors [58] and from the surface diffusion of adsorbates in electronic systems [59]. Both the simulation of the system and length for the forward simulation are 10 and $q = 2$. Additionally, soft start [26], a noise estimation and signal smoothing algorithm, is used to improve the denoising. Figure 15 shows that WmSINDy is the superior method overall due to its higher success rate and better performance across the various error metrics. It is also worth noting that WSINDy and WmSINDy are theoretically based on Gaussian noise to find the parameters $m$ and $p$ for the test function $\phi$. Nevertheless, WmSINDy not only performs well in this system but also provides better results than mSINDy

### 4.6 Identification of noise distributions

As explored in [26], once the noise is characterized, it could be useful to examine the statistical characteristics of that noise, particularly its probability distribution pattern. Understanding the distribution of the noise can provide insights into its nature and help inform strategies for dealing with it. Thus, knowing whether the noise follows a normal distribution, uniform distribution, or some other pattern can guide the selection of appropriate filtering techniques or statistical methods for analysis. The Lorenz system with the same initial conditions and hyperparameters as in Section 4.1.1 is used to investigate the performance of the proposed method for different noise types. The library consists of polynomials up to second order for both mSINDy and WmSINDy methods.

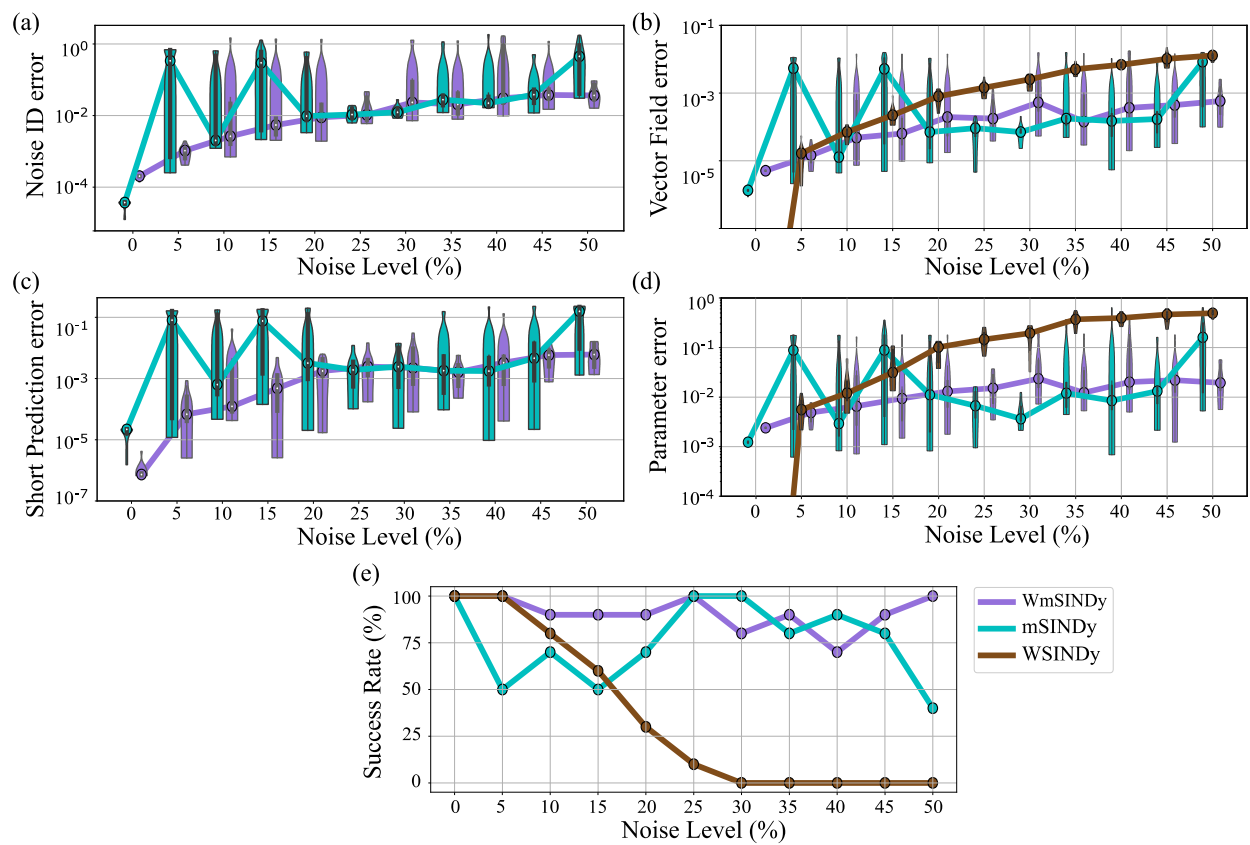

**Fig. 16.** Lorenz system: comparison for different uniform noise levels.

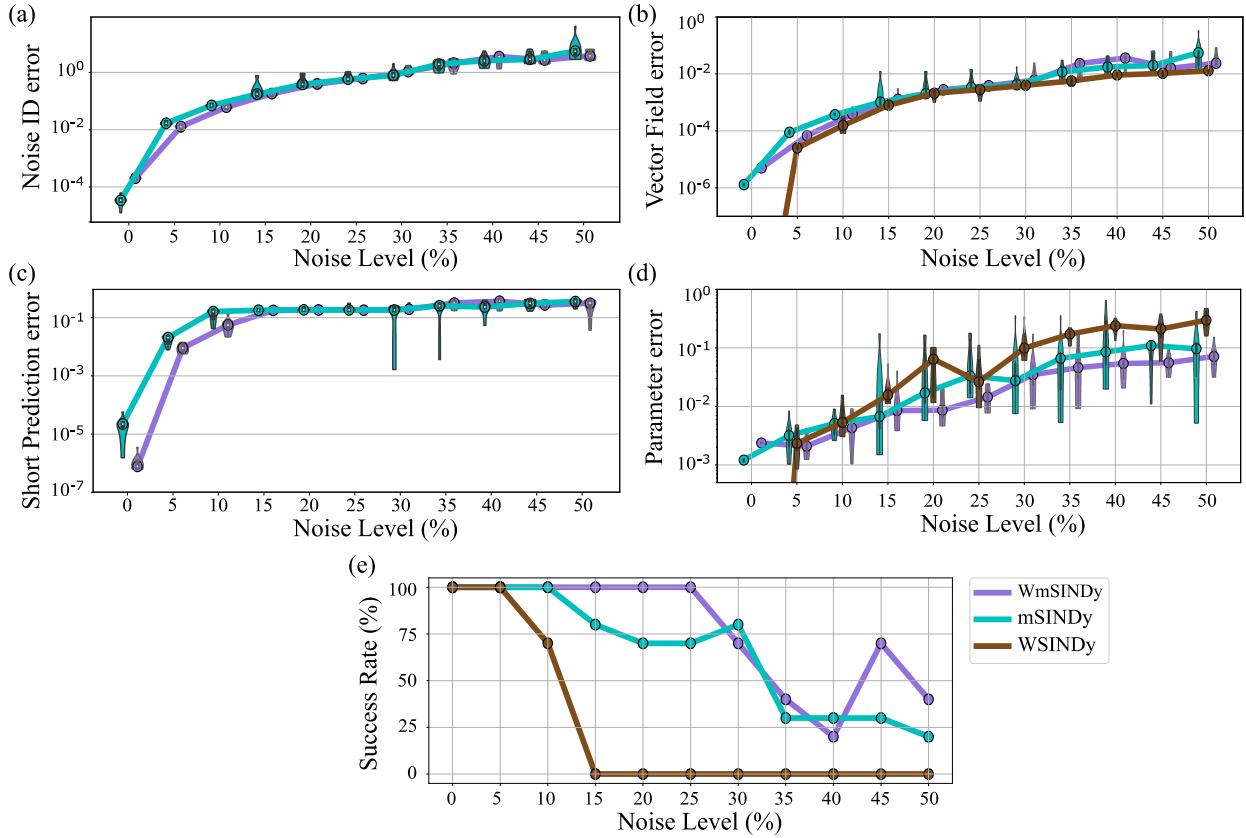

**Fig. 17.** Lorenz system: comparison for different Rayleigh noise levels.

Unlike [26], we use the mean and standard deviations from the normal distribution used to construct Fig. 2 as references. From there, we obtain these similar first two moments for the rest of the distributions, to perform a fair comparison. We employ the **fitter** package [60] to accomplish this task using Python. In Figs. 16-19, we observe the utilization of four distinct noise distributions: uniform, gamma, Dweibull, and Rayleigh, for generating the noisy data. Since each noise distribution has unique properties that model different types of variability or uncertainty, we use these different noise distributions to test the robustness of the WmSINDy method under different types of randomness. We vary the intensity of noise as done in section 4.1.1.

Figure 16 shows the studies made in the presence of uniform noise, a type of noise that is typically observed in images [61]. For the noise error, the WmSINDy generally performs well across all noise levels with lower error values than mSINDy. Also, in WmSINDy there are no abrupt jumps of the medians. Concerning the vector field error, we can see a similar pattern as the noise error. Whereas for the WSINDy, except for the noise-free case, its performance is lower than the other two methods. For the short prediction error, WmSINDy presents better results than mSINDy. For the parameter error, WSINDy starts with a very low noise error but deteriorates as noise increases, with higher errors at higher noise levels. The lowest medians belong to mSINDy; however, it presents considerable variability. WmSINDy shows slightly higher error rates than mSINDy but is more stable than mSINDy. Regarding the success rate, WSINDy starts with a high success rate and drops significantly as noise increases, stabilizing near zero as noise reaches higher levels. mSINDy begins with a high success rate but declines in the low noise level range, then again it reaches high success rates as the noise increases. WmSINDy shows a high success rate, close to 100%, maintaining it throughout the noise levels. In summary, WmSINDy is the most reliable and robust, followed by mSINDy, and finally, WSINDy is the least effective in all conditions, particularly at higher noise levels.

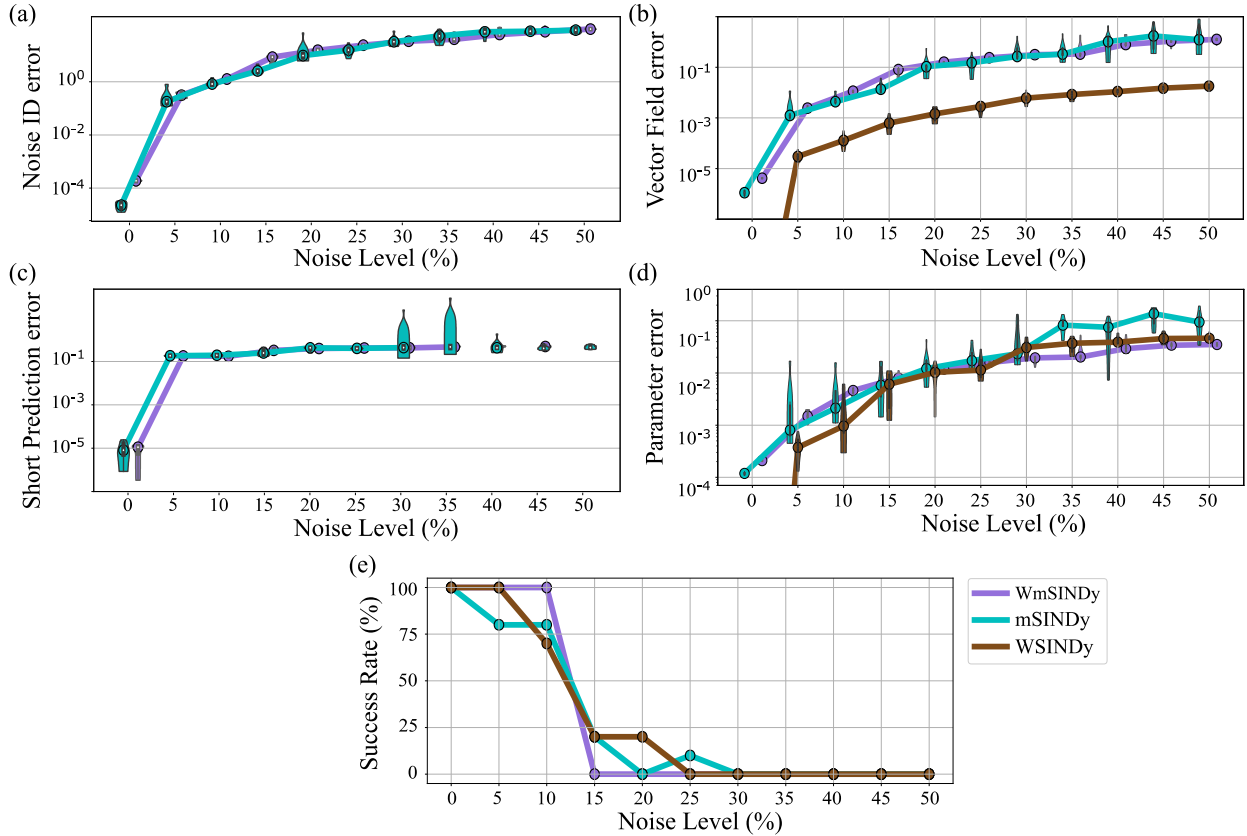

**Fig. 18.** Lorenz system: comparison for different gamma noise levels.

Figure 17 presents the numerical experiments at different levels of Rayleigh noise. This noise is commonly found in Raman distributed temperature sensors[62]. For the noise, vector field, and short prediction errors, the performance of both mSINDy and WmSINDy are similar. While, WSINDy shows the lowest errors at 0 and 5 % levels of noise. As the noise level increases, these errors increase for the three methods. Regarding the parameter error, except for the noise-free case, where WSINDy is the best, WmSINDy shows slightly better performance with respect to the other two methods. For the success rate, WmSINDy provides the best results up to 25% of noise level, then, the behavior is similar for both mSINDy and WmSINDy, with the last one showing more success at high noise levels. WSINDy finds the correct positions of the parameters until a 10% noise level and eventually shows no success beyond this point.

Figure 18 depicts the experiments when using different gamma noise levels. We focus on the noise level of up to 25% of noise, and after that, all methods break. For the noise and vector field errors, mSINDy gives slightly better results than WmSINDy, whereas for the vector field error, WSINDy is considerably better than these two methods The results are similar for the short prediction errors for both mSINDy and WmSINDy. Regarding the parameter error, in the range of noise level 0-10% the best results provide WSINDy, followed by the mSINDy, and finally the WmSINDy, while in the range of noise level 15-25%, all methods perform similarly. As for the success rate, WSINDy is consistently more robust, followed by mSINDy and finally, WmSINDy.

Figure 19 presents the result of using Dweibull noise, a noise commonly observed in ultra-wideband wireless sensor networks[63]. As can be seen, all methods break after 20% of the noise level, a similar behavior was found when using the gamma noise. mSINDy seems to provide slightly better results than the other two methods, followed by WmSINDy, and finally WSINDy.

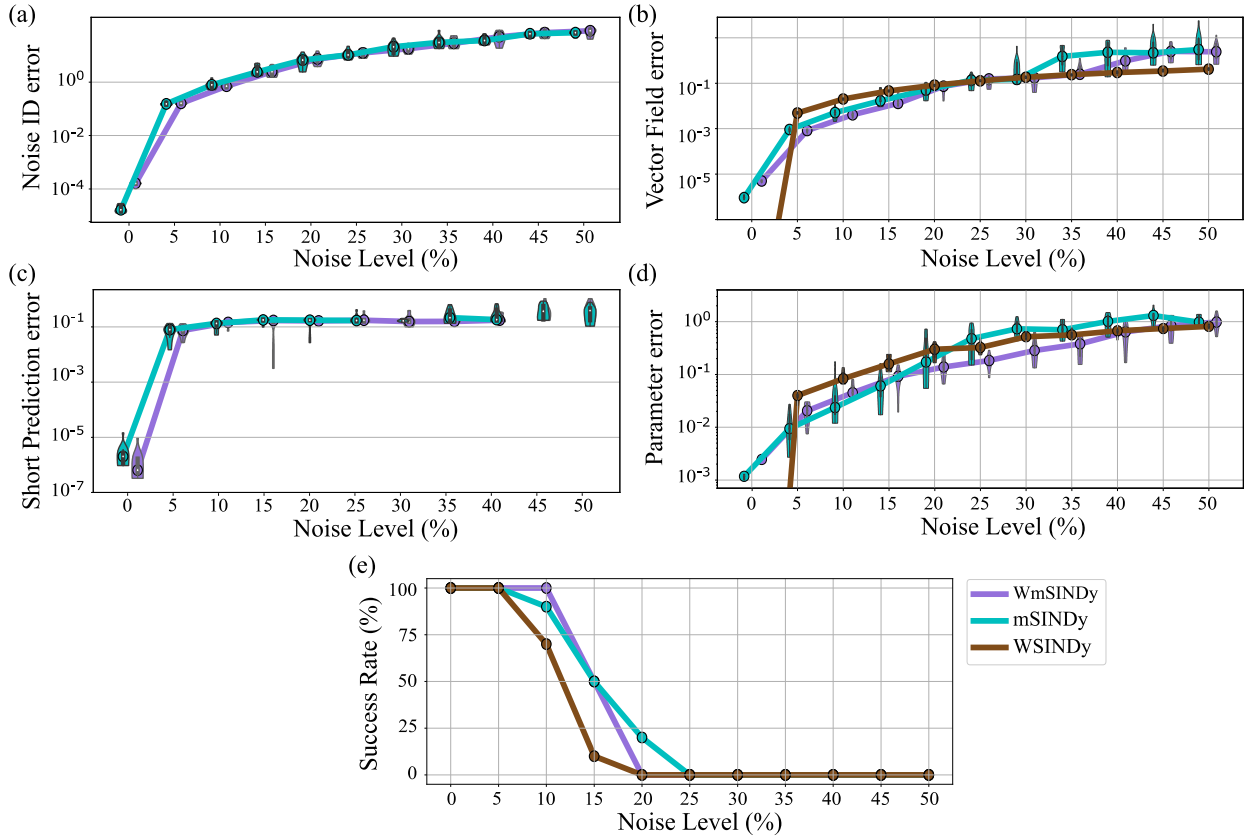

**Fig. 19.** Lorenz system: comparison for different Dweibull noise levels.

When using different noise distributions, such as uniform, gamma, Rayleigh, and Dweibull, with the same mean and standard deviations, the magnitudes of their amplitudes will differ due to their distinct probability density functions. For instance, uniform noise results in consistent perturbations because all values within a given range are equally likely, leading to steady oscillations without significant spikes. With Rayleigh noise, most values are concentrated near the mean, but the distribution has a long tail, so the system may experience moderate fluctuations. Gamma noise, due to its skewness, can produce abrupt jumps in the system, leading to sudden and sharp changes in behavior when noise magnitudes are large. In contrast, Dweibull noise is influenced by its shape parameter, which can yield varying noise characteristics, potentially producing significant extremes that may cause large, sudden deviations in system behavior.

Overall, from Figs. 16-19, at noise-free, WSINDy performs the best followed by mSINDy, and finally WmSINDy. In the low noise level range (up to 15%), WSINDy also provides interesting results. For the rest of the noise intensities, WSINDy shows the lowest variability, but with the highest errors. WmSINDy depicts more variability than WSINDy, but less than mSINDy, since the medians show a monotonic increment as the noise level increases, with no jumps as the behavior of mSINDy. Particularly, the higher success rates of WmSINDy make it more robust than mSINDy and WSINDy to detect the correct dynamics of the system when dealing with different kinds of noise.

The three methods find it particularly challenging to characterize the gamma and Dweibull noise distributions. It is because these noises have a high probability of larger amplitudes, such that they will induce more irregular perturbations and lead to aperiodic behavior, which makes the system identification more difficult.

## 5. Conclusion

We have extended the mSINDy by using the weak formulation from WSINDy to avoid the numerical derivative with the purpose of identifying several dynamical systems embedded in high noise. Overall, the proposed WmSINDy provides better results compared with the baseline methods mSINDy and WSINDy in the different analytical dynamical systems studied. Since WmSINDy is a joint scheme of both mSINDy and WSINDy, we leverage their strengths but also incur some of their weaknesses. For the former, compared to the baseline methods, the proposed algorithm performs well, especially for the systems, such that its configuration exhibits chaotic behavior, at high noise intensities. However, for the systems with limit cycles, the proposed approach performs similarly to mSINDy. Whereas for the latter, we inherit the need to set several hyperparameters. For WmSINDy, we got interesting results using the following values for the hyperparameters from mSINDy. The value of $q$ was set to 1, which also helps to improve the speed of the algorithm. The value of $\lambda$ is correlated with the noise intensity, with higher noise intensity corresponding to higher $\lambda$ values to a certain extent. In our numerical experiments, a value of 0.2 can balance the trade-off between the fit to the data and the number of terms in the model. A small value of $N_{loop}$ is insufficient to characterize both the noise and the dynamics, while values in the range 6 to 8 provide acceptable results. A few data points severely affect the model identification, whereas many data points increase the calculation time. On the other hand, from the WSINDy method, we have focused on only the most influential hyperparameters: $K$ we have set as the length of the number of points, nonetheless, an optimal $K$ must be; $\hat{\tau}$ we set it as -2 as explained in [32], and because the tested systems have no considerable difference between the biggest and smallest values, the scale factor is set to 0. We leave room to find optimal hyperparameters in future work.

The improvement in the accuracy of mSINDy comes at a price in the speed of WmSINDy, because it is slower than both mSINDy and WSINDy. It could be interesting to speed up the proposed technique by using automatic differentiation in JAX M.D. [64,65] or combining it with Rust [66]. Additionally, several approaches have been developed to improve the performance of the Adam optimizer [67,68] Thus, implementing a suitable optimizing technique in the proposed framework might result in better convergence. As was mentioned in [26], combining both mSINDy and WmSINDy with control inputs could open several applications. However, one challenge we encountered was ensuring the optimizer distinguishes between noise and control inputs. Furthermore, it would be valuable to explore an extension of the WmSINDy method in the noise-free cases where the result is the one provided by WSINDy, thus bypassing the need for optimization, as the characterized noise is slightly different from zeros in such cases. Moreover, selecting an appropriate set of candidate terms is crucial for accurate model discovery, particularly when using SINDy and its improvements. Consequently, when using the proposed WmSINDy, one should also pay special attention to capturing the dynamics of the studied systems. Finally, it is worth noting that the performance of both mSINDy and WmSINDy varies depending on the type of noise involved. Specifically, both methods find more challenging scenarios when dealing with non-zero mean noise distributions such as gamma and Rayleigh. Consequently, addressing this limitation in performance represents an important area for future work.